# A RENEWAL THEORY APPROACH TO PERIODIC COPOLYMERS WITH ADSORPTION

By Francesco Caravenna, Giambattista Giacomin[1]
and Lorenzo Zambotti

*Università degli Studi di Padova, Université Paris 7 and Université Paris 6*

We consider a general model of a heterogeneous polymer chain fluctuating in the proximity of an interface between two selective solvents. The heterogeneous character of the model comes from the fact that the *monomer units* interact with the solvents and with the interface according to some *charges* that they carry. The charges repeat themselves along the chain in a periodic fashion. The main question concerning this model is whether the polymer remains tightly close to the interface, a phenomenon called *localization*, or whether there is a marked preference for one of the two solvents, thus yielding a *delocalization* phenomenon.

In this paper, we present an approach that yields sharp estimates for the partition function of the model in all regimes (localized, delocalized and critical). This, in turn, makes possible a precise pathwise description of the polymer measure, obtaining the full *scaling limits* of the model. A key point is the closeness of the polymer measure to suitable *Markov renewal processes*, Markov renewal theory being one of the central mathematical tools of our analysis.

## 1. Introduction and main results.

1.1. *Two motivating models.* Let $S := \{S_n\}_{n=0,1,\ldots}$ be a random walk, $S_0 = 0$ and $S_n = \sum_{j=1}^n X_j$, with i.i.d. symmetric increments taking values in $\{-1, 0, +1\}$. Hence the law of the walk is identified by $p := \mathbf{P}(X_1 = 1) = \mathbf{P}(X_1 = -1)$ and we assume that $p \in (0, 1/2)$. The case $p = 1/2$ can be treated in an analogous way, but requires some notational care because of the periodicity of the walk. We also consider a sequence $\omega := \{\omega_n\}_{n \in \mathbb{N} = \{1,2,\ldots\}}$

Received March 2006; revised March 2006.
[1]Supported in part by ANR project, POLINT BIO.
*AMS 2000 subject classifications.* 60K35, 82B41, 82B44.
*Key words and phrases.* Random walks, renewal theory, Markov renewal theory, scaling limits, polymer models, wetting models.

of real numbers with the property that $\omega_n = \omega_{n+T}$ for some $T \in \mathbb{N}$ and for every $n$, denoting by $T(\omega)$ the minimal value of $T$.

Consider the following two families of modifications of the law of the walk, both indexed by a parameter $N \in \mathbb{N}$:

*Pinning and wetting models.* For $\lambda \geq 0$, consider the probability measure $\mathbf{P}_{N,\omega}$ defined by

$$(1.1) \qquad \frac{\mathrm{d}\mathbf{P}_{N,\omega}}{\mathrm{d}\mathbf{P}}(S) \propto \exp\left(\lambda \sum_{n=1}^{N} \omega_n \mathbf{1}_{\{S_n = 0\}}\right).$$

The walk receives a *pinning reward*, which may be negative or positive, each time it visits the origin. By considering the directed walk viewpoint, that is, $\{(n, S_n)\}_n$, one may interpret this model in terms of a directed linear chain receiving an energetic contribution when it touches an interface. The main question is whether for large $N$ the typical paths of $\mathbf{P}_{N,\omega}$ are attracted or repelled by the interface.

There is an extensive literature on periodic pinning and wetting models, the majority of which is restricted to the $T = 2$ case, for example, [13, 25]; see [15] for further discussion and references.

*Copolymer near a selective interface.* In much the same way, we introduce

$$(1.2) \qquad \frac{\mathrm{d}\mathbf{P}_{N,\omega}}{\mathrm{d}\mathbf{P}}(S) \propto \exp\left(\lambda \sum_{n=1}^{N} \omega_n \operatorname{sign}(S_n)\right),$$

where if $S_n = 0$, we set $\operatorname{sign}(S_n) := \operatorname{sign}(S_{n-1})\mathbf{1}_{\{S_{n-1} \neq 0\}}$. This convention for defining $\operatorname{sign}(0)$, to be kept throughout the paper, simply means that $\operatorname{sign}(S_n) = +1, 0, -1$ according to whether the bond joining $S_{n-1}$ and $S_n$ lies above, on or below the $x$-axis.

Also in this case, we take a directed walk viewpoint and $\mathbf{P}_{N,\omega}$ then may be interpreted as a polymeric chain in which the monomer units, the bonds of the walk, are charged. An interface, the $x$-axis, separates two solvents, say oil above and water below. Positively charged monomers are hydrophobic and negatively charged ones are instead hydrophilic. In this case, one expects competition between three possible scenarios: polymer preferring water, preferring oil or undecided between the two and choosing to fluctuate in the proximity of the interface.

We select [23, 27] from the physical literature on periodic copolymers, keeping in mind, however, plays that periodic copolymer modeling plays a central role in applied chemistry and material science.

1.2. *A general model.* We point out that the models presented in Section 1.1 are particular examples of the polymer measure with Hamiltonian

$$\mathcal{H}_N(S) = \sum_{i=\pm 1} \sum_{n=1}^{N} \omega_n^{(i)} \mathbf{1}_{\{\text{sign}(S_n)=i\}} + \sum_{n=1}^{N} \omega_n^{(0)} \mathbf{1}_{\{S_n=0\}}$$
$$+ \sum_{n=1}^{N} \widetilde{\omega}_n^{(0)} \mathbf{1}_{\{\text{sign}(S_n)=0\}}, \quad (1.3)$$

where $\omega^{(\pm 1)}$, $\omega^{(0)}$ and $\widetilde{\omega}^{(0)}$ are periodic sequences of real numbers. Observe that by our convention on $\text{sign}(0)$, the last term provides an energetic contribution (of pinning/depinning type) to the bonds lying on the interface. We use the shorthand $\omega$ for the four periodic sequences appearing in (1.3) and will use $T = T(\omega)$ to denote the smallest common period of the sequences. We will refer to $\omega$ as to the *charges* of our system.

Besides being a natural model, generalizing and interpolating between pinning and copolymer models, the general model we consider is one which has been considered on several occasions (see, e.g., [28] and references therein).

Starting from the Hamiltonian (1.3), for $a = \text{c}$ (*constrained*) or $a = \text{f}$ (*free*), we introduce the *polymer measure* $\mathbf{P}_{N,\omega}^a$ on $\mathbb{Z}^{\mathbb{N}}$, defined by

$$\frac{d\mathbf{P}_{N,\omega}^a}{d\mathbf{P}}(S) = \frac{\exp(\mathcal{H}_N(S))}{\widetilde{Z}_{N,\omega}^a}(\mathbf{1}_{\{a=\text{f}\}} + \mathbf{1}_{\{a=\text{c}\}}\mathbf{1}_{\{S_N=0\}}), \quad (1.4)$$

where $\widetilde{Z}_{N,\omega}^a := \mathbf{E}[\exp(\mathcal{H}_N)(\mathbf{1}_{\{a=\text{f}\}} + \mathbf{1}_{\{a=\text{c}\}}\mathbf{1}_{\{S_N=0\}})]$ is the *partition function*, that is, the normalization constant. Observe that the polymer measure $\mathbf{P}_{N,\omega}^a$ is invariant under the joint transformation $S \to -S$, $\omega^{(+1)} \to \omega^{(-1)}$, hence, by symmetry, we may (and will) assume that

$$h_\omega := \frac{1}{T(\omega)} \sum_{n=1}^{T(\omega)} (\omega_n^{(+)} - \omega_n^{(-)}) \geq 0. \quad (1.5)$$

We also set $\mathbb{S} := \mathbb{Z}/(T\mathbb{Z})$ and for $\beta \in \mathbb{S}$, we equivalently write $[n] = \beta$ or $n \in \beta$. Notice that the charges $\omega_n$ are functions of $[n]$, so that we can write $\omega_{[n]} := \omega_n$.

1.3. *The free energy viewpoint.* The standard statistical mechanics approach naturally leads to a consideration of the free energy of the model, that is, the limit as $N \to \infty$ of $(1/N) \log \widetilde{Z}_{N,\omega}^a$. It is, however, practical to observe that we can add to the Hamiltonian $\mathcal{H}_N$ a term which is constant with respect to $S$ without changing the polymer measure. Namely, if we set

$$\mathcal{H}_N'(S) := \mathcal{H}_N(S) - \sum_{n=1}^{N} \omega_n^{(+1)},$$

which amounts to sending $\omega_n^{(+1)} \to 0$, $\omega_n^{(-1)} \to (\omega_n^{(-1)} - \omega_n^{(+1)})$ and $\widetilde{\omega}_n^{(0)} \to (\widetilde{\omega}_n^{(0)} - \omega_n^{(+1)})$, then we can write

$$(1.6) \qquad \frac{d\mathbf{P}_{N,\omega}^a}{d\mathbf{P}}(S) = \frac{\exp(\mathcal{H}_N'(S))}{Z_{N,\omega}^a}(\mathbf{1}_{\{a=f\}} + \mathbf{1}_{\{a=c\}}\mathbf{1}_{\{S_N=0\}}),$$

where $Z_{N,\omega}^a$ is a new partition function given by

$$(1.7) \quad Z_{N,\omega}^a = \mathbf{E}[\exp(\mathcal{H}_N')(\mathbf{1}_{\{a=f\}} + \mathbf{1}_{\{a=c\}}\mathbf{1}_{\{S_N=0\}})] = \widetilde{Z}_{N,\omega}^a \cdot e^{-\sum_{n=1}^N \omega_n^{(+1)}}.$$

At this point, we define the free energy:

$$(1.8) \qquad \mathrm{F}_\omega := \lim_{N \to \infty} \frac{1}{N} \log Z_{N,\omega}^c.$$

A proof of the existence of such a limit involves standard superadditive arguments, as well as the fact that the superscript c could be changed to f without changing the result (see, e.g., [15], but a complete proof, without the use of super-additivity, is given below).

The principle that the free energy contains the crucial information on the large $N$ behavior of the system is certainly not violated in this context. In order to clarify this point, let us first observe that $\mathrm{F}_\omega \geq 0$ for every $\omega$. This follows by observing that the energetic contribution to the trajectories that stay positive and come back to zero for the first time at epoch $N$ is just $\omega_N^{(0)}$, hence, by (1.7),

$$(1.9) \qquad \begin{aligned} \frac{1}{N} \log Z_{N,\omega}^c &\geq \frac{\omega_N^{(0)}}{N} + \frac{1}{N} \log \mathbf{P}(S_n > 0, n = 1, \ldots, N-1, S_N = 0) \\ &\stackrel{N \to \infty}{\longrightarrow} 0, \end{aligned}$$

where we have simply used the fact that the distribution of the first return to zero of $S$ is subexponential [see (2.2) for a much sharper estimate]. This suggests a natural dichotomy and, inspired by (1.9), we give the following definition.

DEFINITION 1.1. The polymer chain defined by (1.4) is said to be

- *localized* if $\mathrm{F}_\omega > 0$;
- *delocalized* if $\mathrm{F}_\omega = 0$.

A priori, one is certainly not totally convinced by such a definition. Localization, as well as delocalization, should be given in terms of path properties of the process: it is quite clear that the energy $\mathcal{H}_N'(S)$ of trajectories $S$ which do not come back very often (i.e., not in a positively recurrent fashion) to the interface will either be negative or $o(N)$, but this is far from being a convincing statement of localization. An analogous observation can be made for delocalized polymer chains.

Nonetheless, with a few exceptions, much of the literature focuses on free energy estimates. For example, in [5], one can find the analysis of the free energy of a subset of the class of models we are considering here, and in Section 1.7 of the same work, it is argued that some (weak) path statements of localization and delocalization can be extracted from the free energy. We will return to a review of the existing literature after we have stated our main results, but at this point, we anticipate that our program is going well beyond free energy estimates.

One of the main results in [5] is a formula for $\textsc{f}_\omega$, obtained via large deviations arguments. We will not give the precise expression now (the reader can find it in Section 3.2 below), but we point out that this formula is proved here using arguments that are more elementary and, at the same time, these arguments yield much stronger estimates. More precisely, there exists a positive parameter $\delta^\omega$, given explicitly and analyzed in detail in Section 2.1, that determines the *precise asymptotic behavior* of the partition function (the link between $\delta_\omega$ and $\textsc{f}_\omega$ will be immediately after the statement).

THEOREM 1.2 (Sharp asymptotic estimates). *Fix $\eta \in \mathbb{S}$ and consider the asymptotic behavior of $Z^{\text{c}}_{N,\omega}$ as $N \to \infty$ along $[N] = \eta$. Then:*

(1) *if $\delta^\omega < 1$ then $Z^{\text{c}}_{N,\omega} \sim C^<_{\omega,\eta}/N^{3/2}$;*
(2) *if $\delta^\omega = 1$ then $Z^{\text{c}}_{N,\omega} \sim C^=_{\omega,\eta}/N^{1/2}$;*
(3) *if $\delta^\omega > 1$ then $\textsc{f}_\omega > 0$ and $Z^{\text{c}}_{N,\omega} \sim C^>_{\omega,\eta} \exp(\textsc{f}_\omega N)$,*

*where the quantities $\textsc{f}_\omega$, $C^>_{\omega,\eta}$, $C^<_{\omega,\eta}$ and $C^=_{\omega,\eta}$ are given explicitly in Section 3.*

Of course, by $a_N \sim b_N$ we mean $a_N/b_N \to 1$. We note that Theorem 1.2 implies the existence of the limit in (1.8) and that $\textsc{f}_\omega = 0$ exactly when $\delta^\omega \leq 1$, but we stress that in our arguments, we do not rely on (1.8) to define $\textsc{f}_\omega$. We also point out that analogous asymptotic estimates can be obtained for the free partition function; see Proposition 3.2.

It is rather natural to think that from such precise estimates one can extract detailed information on the limit behaviors of the system. This is correct. Notably, we can consider:

(1) infinite volume limits, that is, weak limits of $\mathbf{P}^a_{N,\omega}$ as a measure on $\mathbb{R}^\mathbb{N}$;
(2) scaling limits, that is, limits in law of the process $S$, suitably rescaled, under $\mathbf{P}^a_{N,\omega}$.

Here, we will focus only on (2); the case (1) is considered in [7].

A word of explanation is in order concerning the fact that there appear to be two types of delocalized polymer chains: those with $\delta^\omega = 1$ and those with $\delta^\omega < 1$. As we will see, these two cases exhibit substantially different path behavior (even if both display distinctive features of delocalized paths, notably a vanishing density of visits at the interface). As will be made clear, in the case $\delta^\omega < 1$, the system is *strictly* delocalized in the sense that a small perturbation in the charges leaves $\delta^\omega < 1$ (as a matter of fact, for charges of a fixed period, the mapping $\omega \mapsto \delta^\omega$ is continuous), while $\delta^\omega$ is rather a borderline, or *critical*, case.

1.4. *The scaling limits.* The main results of this paper concern the *diffusive rescaling* of the polymer measure $\mathbf{P}^a_{N,\omega}$. More precisely, let us define the map $X^N : \mathbb{R}^N \mapsto C([0,1])$:

$$X^N_t(x) = \frac{x_{\lfloor Nt \rfloor}}{\sigma N^{1/2}} + (Nt - \lfloor Nt \rfloor)\frac{x_{\lfloor Nt \rfloor+1} - x_{\lfloor Nt \rfloor}}{\sigma N^{1/2}}, \qquad t \in [0,1],$$

where $\lfloor \cdot \rfloor$ denotes the integer part of $\cdot$ and $\sigma^2 := 2p$ is the variance of $X_1$ under the original random walk measure $\mathbf{P}$. Notice that $X^N_t(x)$ is nothing but the linear interpolation of $\{x_{\lfloor Nt \rfloor}/(\sigma \sqrt{N})\}_{t \in (\mathbb{N}N) \cap [0,1]}$. For $a = \text{f}, \text{c}$ we set

$$Q^a_{N,\omega} := \mathbf{P}^a_{N,\omega} \circ (X^N)^{-1}.$$

Then $Q^a_{N,\omega}$ is a measure on $C([0,1])$, the space of continuous real-valued functions defined on the interval $[0,1]$, and we want to study the behavior as $N \to \infty$ of this sequence of measures.

We start by fixing notation for the following standard processes:

- the Brownian motion $\{B_\tau\}_{\tau \in [0,1]}$;
- the Brownian bridge $\{\beta_\tau\}_{\tau \in [0,1]}$ between 0 and 0;
- the Brownian motion *conditioned to stay nonnegative on* $[0,1]$ or, more precisely, the Brownian meander $\{m_\tau\}_{\tau \in [0,1]}$ (cf. [26]) and its modification by a random flip $\{m^{(p)}_\tau\}_{\tau \in [0,1]}$, defined as $m^{(p)} = \sigma m$, where $\mathbb{P}(\sigma = 1) = 1 - \mathbb{P}(\sigma = -1) = p \in [0,1]$ and $(m, \sigma)$ are independent;
- the Brownian bridge *conditioned to stay nonnegative on* $[0,1]$ or, more precisely, the normalized Brownian excursion $\{e_\tau\}_{\tau \in [0,1]}$, also known as the Bessel bridge of dimension 3 between 0 and 0; see [26]. For $p \in [0,1]$, $\{e^{(p)}_\tau\}_{\tau \in [0,1]}$ is the flipped excursion defined in analogy with $m^{(p)}$;
- the skew Brownian motion $\{B^{(p)}_\tau\}_{\tau \in [0,1]}$ and the skew Brownian bridge $\{\beta^{(p)}_\tau\}_{\tau \in [0,1]}$ of parameter $p$ (cf. [26]) the definitions of which are recalled in Remark 1.5 below.

We introduce a final process, labeled by two parameters $p, q \in [0,1]$: consider a random variable $U \mapsto [0,1]$ with the arcsine law $\mathbb{P}(U \leq t) = \frac{2}{\pi} \arcsin \sqrt{t}$ and processes $\beta^{(p)}, m^{(q)}$ as defined above, with $(U, \beta^{(p)}, m^{(q)})$ an independent triple. We then denote by $\{B^{(p,q)}_\tau\}_{\tau \in [0,1]}$ the process defined by

$$B^{(p,q)}_\tau := \begin{cases} \sqrt{U}\beta^{(p)}_{\tau/U}, & \text{if } \tau \leq U, \\ \sqrt{1-U}m^{(q)}_{(\tau-U)/(1-U)}, & \text{if } \tau > U. \end{cases} \quad (1.10)$$

We then have the following theorem, which is the main result of this paper (see Figure 1).

THEOREM 1.3 (Scaling limits). *For every $\eta \in \mathbb{S}$, if $N \to \infty$ along $[N] = \eta$, then the sequence of measures $\{Q^a_{N,\omega}\}$ on $C([0,1])$ converges weakly. More precisely:*

(1) *for $\delta^\omega < 1$ (strictly delocalized regime), $Q^c_{N,\omega}$ converges to the law of $e^{(\mathsf{p}^c_{\omega,\eta})}$ and $Q^f_{N,\omega}$ converges to the law of $m^{(\mathsf{p}^f_{\omega,\eta})}$ for some parameters $\mathsf{p}^a_{\omega,\eta} \in [0,1]$, $a = \mathrm{f}, \mathrm{c}$;*

(2) *for $\delta^\omega = 1$ (critical regime), $Q^c_{N,\omega}$ converges to the law of $\beta^{(\mathsf{p}_\omega)}$ and $Q^f_{N,\omega}$ converges to the law of $B^{(\mathsf{p}_\omega, \mathsf{q}_{\omega,\eta})}$ for some parameters $\mathsf{p}_\omega, \mathsf{q}_{\omega,\eta} \in [0,1]$;*

(3) *for $\delta^\omega > 1$ (localized regime), $Q^a_{N,\omega}$ converges, as $N \to \infty$, to the measure concentrated on the constant function taking the value zero (no need of the restriction $[N] = \eta$).*

The exact values of the parameters $\mathsf{p}^a_{\omega,\eta}$, $\mathsf{p}_\omega$ and $\mathsf{q}_{\omega,\eta}$ are given in (5.5), (5.7), (5.17) and (5.19). See also Remarks 5.3, 5.4, 5.7 and 5.8.

REMARK 1.4. It is natural to wonder why the results for $\delta^\omega \leq 1$ may depend on $[N] \in \mathbb{S}$. First, we stress that only in very particular cases is there effectively a dependence on $\eta$ and we characterize these instances precisely; see Section 2.3. In particular, there is no dependence on $[N]$ for the two motivating models (pinning and copolymer) described in Section 1.1 and, more generally, if $h_\omega > 0$. However, this *dependence on the boundary condition* phenomenon is not a pathology, but rather a sign of the presence of *first-order phase transitions* in this class of models. Nonetheless, the phenomenon is somewhat surprising since the model is one-dimensional. This issue, which is naturally clarified when dealing with the infinite volume limit of the model, is treated in [7].

REMARK 1.5 (*Skew Brownian motion*). We recall that $B^{(p)}$ (resp. $\beta^{(p)}$) is a process such that $|B^{(p)}| = |B|$ (resp. $|\beta^{(p)}| = |\beta|$) in distribution, but in which the sign of each excursion is chosen to be $+1$ (resp. $-1$) with probability $p$ (resp. $1-p$) instead of $1/2$. Observe that for $p = 1$, we have $B^{(1)} = |B|$, $\beta^{(1)} = |\beta|$, $m^{(1)} = m$ and $e^{(1)} = e$ in distribution. Moreover, $B^{(1/2)} = B$ and $\beta^{(1/2)} = \beta$ in distribution. Notice also that the process $B^{(p,q)}$ differs from the $p$-skew Brownian motion $B^{(p)}$ only for the last excursion in $[0,1]$, whose sign is $+1$ with probability $q$ instead of $p$.

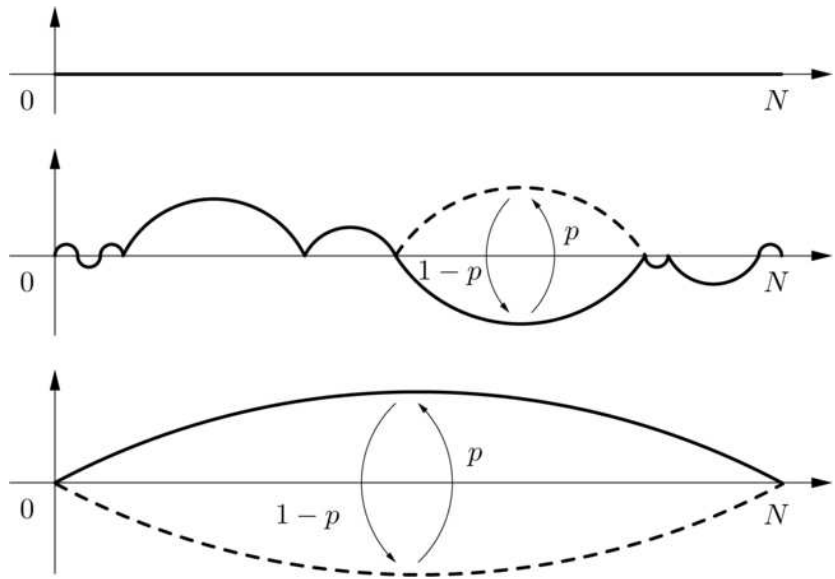

FIG. 1.  *A schematic view of the scaling limits for the constrained endpoint case. While in the localized regime (top image), on a large scale, the polymer cannot be distinguished from the interface, in the strictly delocalized regime (bottom image) the visits to the interface are few and are all close to the endpoints (the sign of the excursion is obtained by flipping one biased coin). In between, there is the critical case: the zeros of the limiting process coincide with the zeros of a Brownian bridge, as found for the homogeneous wetting case [6, 9, 17], but this time, the signs of the excursions vary and they are chosen by flipping independent biased coin. Of course, this suggests that the trajectories in the localized case should be analyzed without rescaling (this is done in [7]).*

1.5. *Motivations and a survey of the literature.* From an applied viewpoint, the interest in periodic models of the type we consider appears to be at least twofold:

(1) On one hand, periodic models are often (e.g., [13, 23]) motivated as caricatures of the *quenched disordered* models, like those in which the charges are a typical realization of a sequence of independent random variables (e.g., [1, 4, 15, 28] and references therein). In this respect, periodic models may be viewed as *weakly inhomogeneous*, and the approximation of *strongly inhomogeneous* quenched models by periodic ones, in the limit of large period, gives rise to very interesting and challenging questions. We believe that if the precise description of the periodic case that we have obtained in this work highlights the limitations of periodic modeling for strongly inhomogeneous systems (cf., in particular, the anomalous decay of quenched partition functions along subsequences pointed out in [16], Section 4 and our Theorem 1.2), it is at the same time an essential step toward understanding the large period limit, and the method we use in this paper may allow a generalization that yields information on this limit.
(2) On the other hand, as mentioned above, periodic models are absolutely natural and of direct relevance in applications, for example, when dealing with *molecularly engineered* polymers; [24, 27] provide a sample of the theoretical physics literature, but the applied literature is extremely vast.

From a mathematical standpoint, our work may be viewed as a further step in the direction of:

(a) extending to the periodic setting precise path estimates obtained for homogeneous models;
(b) clarifying the link between the free energy characterization and the path characterization of the different regimes.

With reference to (a), we point out the novelty with respect to the work on homogeneous models [6, 9, 17, 23]. Although the basic role of renewal theory techniques in obtaining the crucial estimates has already been emphasized in [6, 9], we stress that the underlying key renewal processes that appear in our inhomogeneous context are not standard renewals, but rather *Markov renewal processes* (cf. [2]). Understanding the algebraic structure leading to this type of renewal is one of the central points of our work; see Section 3.1.

We also point out that the Markov renewal processes appearing in the critical regime have step distributions with *infinite mean*. Even for ordinary renewal process, the exact asymptotic behavior of the Green function in the infinite mean case has been a long-standing problem (cf. [14] and [20]) solved only recently by Doney in [10]. The extension of this result to the framework of Markov renewal theory (that we consider here in the case of a finite-state modulating chain) presents some additional problems (see Remark 3.1 and Appendix A) and, to our knowledge, has not been considered in the literature. In Section 5, we also give an extension to our Markov renewal setting of the beautiful theory of convergence of regenerative sets developed in [12].

A final observation is that, as in [6], the estimates we obtain here are really sharp in all regimes and our method goes well beyond the case of random walks with jumps $\pm 1$ and $0$, to which we restrict our attention for the sake of conciseness.

With reference to (b), we point out that in the models we consider, there is a variety of delocalized path behaviors which are not captured by the free energy. This is suggestive also in view of progressing in the understanding of the delocalized phase in the quenched models [16].

1.6. *Outline of the paper.* The remainder of the paper is organized as follows.

- In Section 2, we define the basic parameter $\delta^\omega$ and analyze the dependence of our results on the boundary condition $[N] = \eta$.
- In Section 3, we clarify the connection with Markov renewal theory and obtain the asymptotic behavior of $Z^c_{N,\omega}$ and $Z^f_{N,\omega}$, proving Theorem 1.2.
- In Section 4, we present a basic splitting of the polymer measure into zero level set and excursions and discuss the importance of the partition function.
- In Section 5, we compute the scaling limits of $\mathbf{P}^a_{N,\omega}$, proving Theorem 1.3.
- Finally, the appendices contain the proofs of some technical results.

## 2. A closer look at the main results.

2.1. *The order parameter $\delta^\omega$.* A remarkable feature of our results (see Theorem 1.2 and Theorem 1.3) is the fact that the properties of the polymer measure are essentially encoded in a single parameter $\delta^\omega$ that can be regarded as the *order parameter* of our models. This subsection is devoted to defining this parameter, but we first need some preliminary notation.

We start with the law of the first return to zero of the original walk:

$$(2.1) \qquad \tau_1 := \inf\{n > 0 : S_n = 0\}, \qquad K(n) := \mathbf{P}(\tau_1 = n).$$

It is a classical result [11], Chapter XII.7 that

$$(2.2) \qquad \exists \lim_{n \to \infty} n^{3/2} K(n) =: c_K \in (0, \infty).$$

The key observation is that by the $T$-periodicity of the charges $\omega$ and by the definition (1.5) of $h_\omega$, we can define an $\mathbb{S} \times \mathbb{S}$ matrix $\Sigma_{\alpha,\beta}$ by means of the following relation:

$$(2.3) \qquad \sum_{n=n_1+1}^{n_2} (\omega_n^{(-1)} - \omega_n^{(+1)}) = -(n_2 - n_1)h_\omega + \Sigma_{[n_1],[n_2]}.$$

We have thus decomposed the above sum into a drift term and a fluctuating term, where the latter has the remarkable property of depending on $n_1$ and $n_2$ only through their equivalence classes $[n_1]$ and $[n_2]$ in $\mathbb{S}$. We can now define three basic objects:

- for $\alpha, \beta \in \mathbb{S}$ and $\ell \in \mathbb{N}$, we set

$$(2.4) \qquad \Phi^\omega_{\alpha,\beta}(\ell) := \begin{cases} \omega_\beta^{(0)} + (\widetilde{\omega}_\beta^{(0)} - \omega_\beta^{(+1)}), \\ \qquad \text{if } \ell = 1, \ell \in \beta - \alpha, \\ \omega_\beta^{(0)} + \log[\frac{1}{2}(1 + \exp(-\ell h_\omega + \Sigma_{\alpha,\beta}))], \\ \qquad \text{if } \ell > 1, \ell \in \beta - \alpha, \\ 0, \qquad \text{otherwise;} \end{cases}$$

- for $x \in \mathbb{N}$, we introduce the $\mathbb{S} \times \mathbb{S}$ matrix $M^\omega_{\alpha,\beta}(x)$ defined by

$$M^\omega_{\alpha,\beta}(x) := e^{\Phi^\omega_{\alpha,\beta}(x)} K(x) \mathbf{1}_{(x \in \beta - \alpha)}; \tag{2.5}$$

- summing the entries of $M^\omega$ over $x$, we obtain an $\mathbb{S} \times \mathbb{S}$ matrix that we call $B^\omega$:

$$B^\omega_{\alpha,\beta} := \sum_{x \in \mathbb{N}} M^\omega_{\alpha,\beta}(x). \tag{2.6}$$

The meaning of these quantities will become clear in the next subsection. For now we stress that they are *explicit functions* of the charges $\omega$ and of the law of the underlying random walk (to ease notation, the $\omega$-dependence of these quantities will often be dropped).

Observe that $B_{\alpha,\beta}$ is a finite-dimensional matrix with positive entries, hence the Perron–Frobenius theorem (see, e.g., [2]) implies that $B_{\alpha,\beta}$ has a unique real positive eigenvalue (called the *Perron–Frobenius eigenvalue*) with the property that it is a simple root of the characteristic polynomial and that it coincides with the spectral radius of the matrix. This is exactly our order parameter:

$$\delta^\omega := \text{Perron–Frobenius eigenvalue of } B^\omega. \tag{2.7}$$

2.2. *A random walk excursions viewpoint.* In this subsection, we are going to see that the quantities defined in (2.4) and (2.5) emerge in a natural way from the algebraic structure of the constrained partition function $Z^c_{N,\omega}$. Let us reconsider our Hamiltonian (1.3): its specificity is due to the fact that it can be decomposed in an efficient way by considering the return times to the origin of $S$. More precisely, we define

$$\tau_0 = 0, \qquad \tau_{j+1} = \inf\{n > \tau_j : S_n = 0\}$$

for $j \in \mathbb{N}$ and we set $\iota_N = \sup\{k : \tau_k \leq N\}$. We also set $T_j = \tau_j - \tau_{j-1}$ and of course $\{T_j\}_{j=1,2,\ldots}$ is, under $\mathbf{P}$, an i.i.d. sequence. By conditioning on $\tau$ and integrating on the up–down symmetry of the random walk excursions, one easily obtains the following expression for the constrained partition function:

$$\begin{aligned}
Z^c_{N,\omega} &= \mathbf{E}\left[\prod_{j=1}^{\iota_N} \exp(\Phi^\omega_{[\tau_{j-1}],[\tau_j]}(\tau_j - \tau_{j-1})); \tau_{\iota_N} = N\right] \\
&= \sum_{k=1}^N \sum_{\substack{t_0,\ldots,t_k \in \mathbb{N} \\ 0 =: t_0 < t_1 < \cdots < t_k := N}} \prod_{j=1}^k M^\omega_{[t_{j-1}],[t_j]}(t_j - t_{j-1})
\end{aligned} \tag{2.8}$$

where we have used the quantities introduced in (2.4) and (2.5). This formula shows in particular that the partition function $Z^c_{N,\omega}$ is a function of the entries of $M^\omega$.

We stress that the algebraic form of (2.8) is of crucial importance. It will be analyzed in detail and exploited in Section 3 and will be the key to the proof of Theorem 1.2.

2.3. *The regime* $\omega \in \mathcal{P}$. In this subsection, we look more closely at the dependence of our main result (Theorem 1.3) on the boundary condition $[N] = \eta$. It is convenient to introduce the subset $\mathcal{P}$ of charges defined by

(2.9) $\qquad \mathcal{P} := \{\omega : \delta^\omega \leq 1, h_\omega = 0, \exists \alpha, \beta : \Sigma_{\alpha,\beta} \neq 0\},$

where we recall that $h_\omega$ and $\Sigma_{\alpha,\beta}$ have been defined respectively in (1.5) and (2.3).

The basic observation is that if $\omega \notin \mathcal{P}$, the constants $\mathtt{p}^{\mathrm{c}}_{\omega,\eta}$, $\mathtt{p}^{\mathrm{f}}_{\omega,\eta}$, $\mathtt{p}_\omega$ and $\mathtt{q}_{\omega,\eta}$ actually have no dependence on $\eta$ and all take the same value, namely 1 if $h_\omega > 0$ and $1/2$ if $h_\omega = 0$ (see Remarks 5.3, 5.4, 5.7 and 5.8). The results in Theorem 1.3 for $\delta^\omega \leq 1$ may then be strengthened in the following way:

PROPOSITION 2.1. *If $\omega \notin \mathcal{P}$, then the sequence of measures $\{Q^a_{N,\omega}\}$ on $C([0,1])$ converges weakly as $N \to \infty$. In particular, setting $\mathtt{p}_\omega := 1$ if $h_\omega > 0$ and $\mathtt{p}_\omega := \frac{1}{2}$ if $h_\omega = 0$:*

(1) *for $\delta^\omega < 1$ (strictly delocalized regime), $Q^{\mathrm{f}}_{N,\omega}$ converges to the law of $m^{(\mathtt{p}_\omega)}$ and $Q^{\mathrm{c}}_{N,\omega}$ converges to the law of $e^{(\mathtt{p}_\omega)}$;*

(2) *for $\delta^\omega = 1$ (critical regime), $Q^{\mathrm{f}}_{N,\omega}$ converges to the law of $B^{(\mathtt{p}_\omega)}$ and $Q^{\mathrm{c}}_{N,\omega}$ converges to the law of $\beta^{(\mathtt{p}_\omega)}$.*

This stronger form of the scaling limits holds, in particular, for the two motivating models of Section 1.1, the pinning and the copolymer models, for which $\omega$ *never* belongs to $\mathcal{P}$. This is clear for the pinning case, where, by definition, $\Sigma \equiv 0$, while the copolymer model with $h_\omega = 0$ always has $\delta^\omega > 1$, as we will prove in Appendix B. However, we stress that there do exist charges $\omega$ (necessarily belonging to $\mathcal{P}$) for which there is indeed a dependence on $[N] = \eta$ in the delocalized and critical scaling limits. This interesting phenomenon may be understood in statistical mechanics terms and is analyzed in detail in [7].

**3. Sharp asymptotic behavior for the partition function.** In this section, we are going to derive the precise asymptotic behavior of $Z^{\mathrm{c}}_{N,\omega}$ and $Z^{\mathrm{f}}_{N,\omega}$, in particular, proving Theorem 1.2. The key observation is that the study of the partition function for the models we are considering can be seen as part of the framework of the theory of *Markov renewal processes*; see [2], Chapter VII.4.

3.1. *A link with Markov renewal theory.* The starting point of our analysis is equation (2.8). Let us call a function $\mathbb{N} \times \mathbb{S} \times \mathbb{S} \ni (x, \alpha, \beta) \mapsto F_{\alpha,\beta}(x) \geq 0$ a *kernel*. For fixed $x \in \mathbb{N}$, $F_{\cdot,\cdot}(x)$ is an $\mathbb{S} \times \mathbb{S}$ matrix with nonnegative entries. Given two kernels $F$ and $G$, we define their convolution $F * G$ as the kernel defined by

$$(3.1) \quad (F * G)_{\alpha,\beta}(x) := \sum_{y \in \mathbb{N}} \sum_{\gamma \in \mathbb{S}} F_{\alpha,\gamma}(y) G_{\gamma,\beta}(x-y) = \sum_{y \in \mathbb{N}} [F(y) \cdot G(x-y)]_{\alpha,\beta},$$

where $\cdot$ denotes matrix product. Then, since, by construction, $M_{\alpha,\beta}(x) = 0$ if $[x] \neq \beta - \alpha$, we can write (2.8) in the following way:

$$(3.2) \quad \begin{aligned} Z^{\mathrm{c}}_{N,\omega} &= \sum_{k=1}^{N} \sum_{\substack{t_1,\ldots,t_k \in \mathbb{N} \\ 0 < t_1 < \cdots < t_k := N}} [M(t_1) \cdot \ldots \cdot M(N - t_{k-1})]_{0,[N]} \\ &= \sum_{k=1}^{\infty} [M^{*k}]_{[0],[N]}(N), \end{aligned}$$

where $F^{*n}$ denotes the $n$-fold convolution of a kernel $F$ with itself (the $n = 0$ case is, by definition, $[F^{*0}]_{\alpha,\beta}(x) := \mathbf{1}_{(\beta = \alpha)} \mathbf{1}_{(x=0)}$). In view of (3.2), we introduce the kernel

$$(3.3) \quad \mathcal{Z}_{\alpha,\beta}(n) := \sum_{k=1}^{\infty} [M^{*k}]_{\alpha,\beta}(n)$$

so that $Z^{\mathrm{c}}_{N,\omega} = \mathcal{Z}_{[0],[N]}(N)$ and, more generally, $Z^{\mathrm{c}}_{N-k,\theta_k \omega} = \mathcal{Z}_{[k],[N]}(N-k)$, $k \leq N$, where we have introduced the shift operator for $k \in \mathbb{N}$,

$$\theta_k : \mathbb{R}^{\mathbb{S}} \mapsto \mathbb{R}^{\mathbb{S}}, \qquad (\theta_k \zeta)_\beta := \zeta_{[k]+\beta}, \qquad \beta \in \mathbb{S}.$$

Our goal is to determine the asymptotic behavior as $N \to \infty$ of the kernel $\mathcal{Z}_{\alpha,\beta}(N)$ and hence of the partition function $Z^{\mathrm{c}}_{N,\omega}$. To this end, we introduce an important transformation of the kernel $M$ that exploits the algebraic structure of (3.3): we suppose that $\delta^\omega \geq 1$ (the case $\delta^\omega < 1$ requires a different procedure) and set for $b \geq 0$ (cf. [2], Theorem 4.6)

$$A^b_{\alpha,\beta}(x) := M_{\alpha,\beta}(x) e^{-bx}.$$

Let us denote by $\Delta(b)$ the Perron–Frobenius eigenvalue of the matrix $\sum_x A^b_{\alpha,\beta}(x)$. As the entries of this matrix are analytic and nonincreasing functions of $b$, $\Delta(b)$ is analytic and nonincreasing too, hence strictly decreasing, because $\Delta(0) = \delta^\omega \geq 1$ and $\Delta(\infty) = 0$. Therefore, there exists a single value $\mathrm{F}_\omega \geq 0$ such that $\Delta(\mathrm{F}_\omega) = 1$ and we denote by $\{\zeta_\alpha\}_\alpha$, $\{\xi_\alpha\}_\alpha$ the Perron–Frobenius left and right eigenvectors of $\sum_x A^{\mathrm{F}_\omega}_{\alpha,\beta}(x)$, chosen to have (strictly) positive components and normalized in such a way that $\sum_\alpha \zeta_\alpha \xi_\alpha = 1$ (the remaining degree of freedom in the normalization is immaterial). We then set

$$\Gamma_{\alpha,\beta}(x) := M_{\alpha,\beta}(x) e^{-\mathrm{F}_\omega x} \frac{\xi_\beta}{\xi_\alpha} \tag{3.4}$$

and observe that we can rewrite (3.2) as

$$\mathcal{Z}_{\alpha,\beta}(n) = \exp(\mathrm{F}_\omega n) \frac{\xi_\alpha}{\xi_\beta} U_{\alpha,\beta}(n) \qquad \text{where } U_{\alpha,\beta}(n) := \sum_{k=1}^{\infty} [\Gamma^{*k}]_{\alpha,\beta}(n). \tag{3.5}$$

The kernel $U_{\alpha,\beta}(n)$ has a basic probabilistic interpretation that we now describe. Notice first that by construction, we have $\sum_{\beta,x} \Gamma_{\alpha,\beta}(x) = 1$, that is, $\Gamma$ is a *semi-Markov kernel* (cf. [2]). We can then define a Markov chain $\{(J_k, T_k)\}$ on $\mathbb{S} \times \mathbb{N}$ by

$$\mathbb{P}[(J_{k+1}, T_{k+1}) = (\beta, x)|(J_k, T_k) = (\alpha, y)] = \Gamma_{\alpha,\beta}(x) \tag{3.6}$$

and we denote by $\mathbb{P}_\alpha$ the law of $\{(J_k, T_k)\}$ with starting point $J_0 = \alpha$ (the value of $T_0$ plays no role). The probabilistic meaning of $U_{\alpha,\beta}(x)$ is then given by

$$U_{\alpha,\beta}(n) = \sum_{k=1}^{\infty} \mathbb{P}_\alpha(T_1 + \cdots + T_k = n, J_k = \beta). \tag{3.7}$$

We point out that the process $\{\tau_k\}_{k \geq 0}$ defined by $\tau_0 := 0$ and $\tau_k := T_1 + \cdots + T_k$ under the law $\mathbb{P}_\alpha$ is what is called a (discrete) *Markov renewal process* (cf. [2]). This provides a generalization of classical renewal processes since the increments $\{T_k\}$ are not i.i.d. but rather are governed by the process $\{J_k\}$ in the way prescribed by (3.6). The process $\{J_k\}$ is called the *modulating chain* and it is indeed a genuine Markov chain on $\mathbb{S}$, with transition kernel $\sum_{x \in \mathbb{N}} \Gamma_{\alpha,\beta}(x)$, while in general, the process $\{T_k\}$ is *not* a Markov chain. One can view $\tau = \{\tau_n\}$ as a (random) subset of $\mathbb{N}$. More generally, it is convenient to introduce the subset

$$\tau^\beta := \bigcup_{k \geq 0 : J_k = \beta} \{\tau_k\}, \qquad \beta \in \mathbb{S}, \tag{3.8}$$

so that equation (3.7) can be rewritten as

$$U_{\alpha,\beta}(n) = \mathbb{P}_\alpha(n \in \tau^\beta). \tag{3.9}$$

This shows that the kernel $U_{\alpha,\beta}(n)$ is really an extension of the Green function of a classical renewal process. In analogy with the classical case, the asymptotic behavior of $U_{\alpha,\beta}(n)$ is sharply linked to the asymptotic behavior of the kernel $\Gamma$, that is, of $M$. To this end, we notice that our setting is a *heavy-tailed* one. More precisely, for every $\alpha, \beta \in \mathbb{S}$, by (2.2), (2.5) and (2.4), we have

$$
\begin{aligned}
(3.10) \quad &\lim_{\substack{x \to \infty \\ [x] = \beta - \alpha}} x^{3/2} M_{\alpha,\beta}(x) \\
&= L_{\alpha,\beta} := \begin{cases} c_K \frac{1}{2}(1 + \exp(\Sigma_{\alpha,\beta})) \exp(\omega_\beta^{(0)}), & \text{if } h_\omega = 0, \\ c_K \frac{1}{2} \exp(\omega_\beta^{(0)}), & \text{if } h_\omega > 0. \end{cases}
\end{aligned}
$$

The rest of this section is devoted to determining the asymptotic behavior of $U_{\alpha,\beta}(n)$ and hence of $\mathcal{Z}_{\alpha,\beta}(n)$, in particular, proving Theorem 1.2. For convenience, we consider the three regimes separately.

3.2. *The localized regime* ($\delta^\omega > 1$). If $\delta^\omega > 1$ then necessarily $\mathrm{F}_\omega > 0$. Notice that $\sum_x \Gamma_{\alpha,\beta}(x) > 0$, so that in particular the modulating chain $\{J_k\}$ is irreducible. The unique invariant measure $\{\nu_\alpha\}_\alpha$ is easily seen to be equal to $\{\zeta_\alpha \xi_\alpha\}_\alpha$.

Let us compute the *mean* $\mu$ of the semi-Markov kernel $\Gamma$:

$$
\mu := \sum_{\alpha,\beta \in \mathbb{S}} \sum_{x \in \mathbb{N}} x \nu_\alpha \Gamma_{\alpha,\beta}(x) = \sum_{\alpha,\beta \in \mathbb{S}} \sum_{x \in \mathbb{N}} x e^{-\mathrm{F}_\omega x} \zeta_\alpha M_{\alpha,\beta}(x) \xi_\beta
$$
$$
= -\frac{\partial \Delta}{\partial b}\bigg|_{b=\mathrm{F}_\omega} \in (0, \infty)
$$

(for the last equality, see, e.g., [5], Lemma 2.1). We can then apply the Markov renewal theorem (cf. [2], Theorem VII.4.3) which in our periodic setting gives

$$
(3.11) \quad \exists \lim_{\substack{x \to \infty \\ [x] = \beta - \alpha}} U_{\alpha,\beta}(x) = T \frac{\nu_\beta}{\mu}.
$$

By (3.5), we then obtain the desired asymptotic behavior:

$$
(3.12) \quad \mathcal{Z}_{\alpha,\beta}(x) \sim \xi_\alpha \zeta_\beta \frac{T}{\mu} \exp(\mathrm{F}_\omega x), \qquad x \to \infty, [x] = \beta - \alpha,
$$

and for $\alpha = [0]$ and $\beta = \eta$, we have proven part (3) of Theorem 1.2, with $C^>_{\omega,\eta} = \xi_0 \zeta_\eta T/\mu$.

3.3. *The critical case* ($\delta^\omega = 1$). In this case, $\mathrm{F}_\omega = 0$ and equation (3.4) reduces to

$$
(3.13) \quad \Gamma_{\alpha,\beta}(x) = M_{\alpha,\beta}(x) \frac{\xi_\beta}{\xi_\alpha}.
$$

The random set $\tau^\beta$ introduced in (3.8) can be written as the union $\tau^\beta = \bigcup_{k \geq 0} \{\tau_k^\beta\}$, where the points $\{\tau_k^\beta\}_{k \geq 0}$ are taken in increasing order, and we set $T_k^\beta := \tau_{k+1}^\beta - \tau_k^\beta$ for $k \geq 0$. Notice that the increments $\{T_k^\beta\}$ correspond to sums of the variables $\{T_i\}$ between the visits of the chain $\{J_k\}$ to the state $\beta$: for instance, we have

$$T_1^\beta = T_{\kappa+1} + \cdots + T_\ell,$$
$$\kappa := \inf\{k \geq 0 : J_k = \beta\}, \ell := \inf\{k > \kappa : J_k = \beta\}.$$

Equation (3.6) then yields that $\{T_k^\beta\}_{k \geq 0}$ is an independent sequence under $\mathbb{P}_\alpha$ and that for $k \geq 1$, the variables $T_k^\beta$ have the same distribution $q^\beta(n) := \mathbb{P}_\alpha(T_1^\beta = n)$ that does not depend on $\alpha$. In general, the variable $T_0^\beta$ has a different law, $q^{(\alpha;\beta)}(n) := \mathbb{P}_\alpha(T_0^\beta = n)$.

These considerations yield the following crucial observation: for fixed $\alpha$ and $\beta$, the process $\{\tau_k^\beta\}_{k \geq 0}$ under $\mathbb{P}_\alpha$ is a (delayed) classical renewal process, with typical step distribution $q^\beta(\cdot)$ and initial step distribution $q^{(\alpha;\beta)}(\cdot)$. By (3.9), $U_{\alpha,\beta}(n)$ is nothing but the Green function (or renewal mass function) of this process; more explicitly, we can write

$$(3.14) \qquad U_{\alpha,\beta}(x) = \left(q^{(\alpha;\beta)} * \sum_{n=0}^\infty (q^\beta)^{*n}\right)(x).$$

Of course, $q^{(\alpha;\beta)}$ plays no role in the asymptotic behavior of $U_{\alpha,\beta}(x)$. The key point is rather the precise asymptotic behavior of $q^\beta(x)$ as $x \to \infty$, $x \in [0]$, which is given by

$$(3.15) \qquad q^\beta(x) \sim \frac{c_\beta}{x^{3/2}} \qquad \text{where } c_\beta := \frac{1}{\zeta_\beta \xi_\beta} \sum_{\alpha,\gamma} \zeta_\alpha L_{\alpha,\gamma} \xi_\gamma > 0,$$

as is proved in detail in Appendix A. The asymptotic behavior of (3.14) then follows by a result of Doney (cf. [10], Theorem B):

$$(3.16) \qquad U_{\alpha,\beta}(x) \sim \frac{T^2}{2\pi c_\beta} \frac{1}{\sqrt{x}}, \qquad x \to \infty, [x] = \beta - \alpha$$

(the factor $T^2$ is due to our periodic setting). Combining equations (3.5), (3.15) and (3.16), we finally get the asymptotic behavior of $\mathcal{Z}_{\alpha,\beta}(x)$:

$$(3.17) \quad \mathcal{Z}_{\alpha,\beta}(x) \sim \frac{T^2}{2\pi} \frac{\xi_\alpha \zeta_\beta}{\sum_{\gamma,\gamma'} \zeta_\gamma L_{\gamma,\gamma'} \xi_{\gamma'}} \frac{1}{\sqrt{x}}, \qquad x \to \infty, [x] = \beta - \alpha.$$

Taking $\alpha = [0]$ and $\beta = \eta$, we have the proof of part (2) of Theorem 1.2.

REMARK 3.1. We point out that formula (3.15) is quite nontrivial. First, the asymptotic behavior $x^{-3/2}$ of the law of the variables $T_1^\beta$ is the same as that of the $T_i$, although $T_1^\beta$ is the sum of a random number of the nonindependent variables $(T_i)$. Second, the computation of the prefactor $c_\beta$ is by no means an obvious task (we stress that the precise value of $c_\beta$ is crucial in the proof of Proposition 5.5 below).

3.4. *The strictly delocalized case* ($\delta^\omega < 1$). We prove that the asymptotic behavior of $\mathcal{Z}_{\alpha,\beta}(x)$ when $\delta^\omega < 1$ is given by

$$\mathcal{Z}_{\alpha,\beta}(x) \sim ([(1-B)^{-1}L(1-B)^{-1}]_{\alpha,\beta})\frac{1}{x^{3/2}},$$
(3.18)
$$x \to \infty, [x] = \beta - \alpha,$$

where the matrices $L$ and $B$ have been defined in (3.10) and (2.6). In particular, taking $\alpha = [0]$ and $\beta = \eta$, (3.18) proves part (1) of Theorem 1.2 with

$$C^<_{\omega,\eta} := [(1-B)^{-1}L(1-B)^{-1}]_{0,\eta}.$$

To start with, it is easily checked by induction that for every $n \in \mathbb{N}$,

(3.19)
$$\sum_{x \in \mathbb{N}}[M^{*n}]_{\alpha,\beta}(x) = [B^n]_{\alpha,\beta}.$$

Next, we claim that, by (3.10), for every $\alpha, \beta \in \mathbb{S}$,

(3.20)
$$\exists \lim_{\substack{x \to \infty \\ [x] = \beta - \alpha}} x^{3/2}[M^{*k}]_{\alpha,\beta}(x) = \sum_{i=0}^{k-1}[B^i \cdot L \cdot B^{(k-1)-i}]_{\alpha,\beta}.$$

We proceed by induction on $k$. The $k = 1$ case is just equation (3.10) and we have that

$$M^{*(n+1)}(x) = \sum_{y=1}^{x/2}(M(y) \cdot M^{*n}(x-y) + M(x-y) \cdot M^{*n}(y))$$

(strictly speaking this formula is true only when $x$ is even, however the odd $x$ case is analogous). By the induction hypothesis, equation (3.20) holds for every $k \leq n$ and, in particular, this implies that $\{x^{3/2}[M^{*k}]_{\alpha,\beta}(x)\}_{x \in \mathbb{N}}$ is a bounded sequence. Therefore, we can apply the dominated convergence theorem and using (3.19), we get

$$\exists \lim_{\substack{x \to \infty \\ [x] = \beta - \alpha}} x^{3/2}[M^{*(n+1)}]_{\alpha,\beta}(x)$$

$$= \sum_\gamma \left(B_{\alpha,\gamma}\sum_{i=0}^{n-1}[B^i \cdot L \cdot B^{(n-1)-i}]_{\gamma,\beta} + L_{\alpha,\gamma}[B^{*n}]_{\gamma,\beta}\right)$$

$$= \sum_{i=0}^{n}[B^i \cdot L \cdot B^{n-i}]_{\alpha,\beta}.$$

Our purpose is to apply the asymptotic result (3.20) to the terms of (3.5) and we need a bound to apply the dominated convergence theorem. We are going to show that

(3.21)
$$x^{3/2}[M^{*k}]_{\alpha,\beta}(x) \leq Ck^3[B^k]_{\alpha,\beta}$$

for some positive constant $C$ and for all $\alpha, \beta \in \mathbb{S}$ and $x, k \in \mathbb{N}$. We again proceed by induction: for the $k = 1$ case, thanks to (3.10), it is possible to find $C$ such that (3.21) holds (this fixes $C$ once and for all). Now, assuming that (3.21) holds for all $k < n$, we show that it also does for $k = n$ (we suppose, for simplicity, that $n = 2m$ is even, the odd $n$ case being analogous). We then have (also assuming that $x$ is even, for simplicity)

$$x^{3/2}[M^{*2m}]_{\alpha,\beta}(x)$$
$$= 2 \sum_{y=1}^{x/2} \sum_{\gamma \in \mathbb{S}} [M^{*m}]_{\alpha,\gamma}(y) x^{3/2} [M^{*m}]_{\gamma,\beta}(x-y)$$
$$\leq 2 \cdot 2^{3/2} C m^3 \sum_{y=1}^{x/2} \sum_{\gamma \in \mathbb{S}} [M^{*m}]_{\alpha,\gamma}(y) [B^m]_{\gamma,\beta} \leq C(2m)^3 [B^{2m}]_{\alpha,\beta},$$

where we have applied (3.19), and (3.21) is proved.

The r.h.s. of (3.21) is summable in $k$ because the matrix $B$ has spectral radius $\delta^\omega < 1$. We can thus apply the dominated convergence theorem to (3.5) using (3.20) and we obtain (3.18) by

$$\lim_{\substack{x \to \infty \\ [x] = \beta - \alpha}} x^{3/2} \mathcal{Z}_{\alpha,\beta}(x) = \sum_{k=1}^{\infty} \sum_{i=0}^{k-1} [B^i \cdot L \cdot B^{(k-1)-i}]_{\alpha,\beta}$$
$$= [(1-B)^{-1} \cdot L \cdot (1-B)^{-1}]_{\alpha,\beta}.$$

This concludes the proof of Theorem 1.3. □

3.5. *The free partition function.* We now want to compute the asymptotic behavior of the free partition function. In particular, we have the following:

PROPOSITION 3.2 (Sharp asymptotic estimates, free case). *As $N \to \infty$, $[N] = \eta$, we have:*

(1) *for $\delta^\omega < 1$ (strictly delocalized regime), $Z_{N,\omega}^{\mathrm{f}} \sim C_{\omega,\eta}^{<,\mathrm{f}}/N^{1/2}$;*
(2) *for $\delta^\omega = 1$ (critical regime), $Z_{N,\omega}^{\mathrm{f}} \sim C_{\omega,\eta}^{=,\mathrm{f}}$;*
(3) *for $\delta^\omega > 1$ (localized regime), $Z_{N,\omega}^{\mathrm{f}} \sim C_{\omega,\eta}^{>,\mathrm{f}} \exp(\mathrm{F}_\omega N)$,*

*where $C_{\omega,\eta}^{>,\mathrm{f}}$, $C_{\omega,\eta}^{<,\mathrm{f}}$ and $C_{\omega,\eta}^{=,\mathrm{f}}$ are explicit positive constants, depending on $\omega$ and $\eta$.*

PROOF. Conditioning on the last zero of $S$ before epoch $N$, we have the useful formula

$$(3.22) \qquad Z_{N,\omega}^{\mathrm{f}} = \sum_{t=0}^{N} Z_{t,\omega}^{\mathrm{c}} P(N-t) \exp(\widetilde{\Phi}_{[t],[N]}(N-t)),$$

where $P(n) := \mathbf{P}(\tau_1 > n) = \sum_{k=n+1}^{\infty} K(k)$ and

(3.23) $\quad \widetilde{\Phi}_{\alpha,\beta}(\ell) := \log[\frac{1}{2}(1 + \exp(-\ell h_\omega + \Sigma_{\alpha,\beta}))]\mathbf{1}_{(\ell>1)}\mathbf{1}_{(\ell\in\beta-\alpha)},$

which differs from $\Phi$ in not having the terms of interaction with the interface [cf. (2.4)].

Since the asymptotic behavior of $P(\ell)\exp(\widetilde{\Phi}_{\alpha,\beta}(\ell))$ will be also needed, we set

$$
\begin{aligned}
\widetilde{L}_{\alpha,\beta} &:= \lim_{\ell\to\infty, \ell\in\beta-\alpha} \sqrt{\ell} P(\ell) e^{\widetilde{\Phi}_{\alpha,\beta}(\ell)} \\
&= \begin{cases} c_K(1 + \exp(\Sigma_{\alpha,\beta})), & \text{if } h_\omega = 0, \\ c_K, & \text{if } h_\omega > 0, \end{cases}
\end{aligned}
\tag{3.24}
$$

as it follows easily from (3.23) and from the fact that $P(\ell) \sim 2c_K/\sqrt{\ell}$ as $\ell \to \infty$. For the rest of the proof, we consider the different regimes separately.

*The strictly delocalized case.* Notice that

$$N^{1/2} Z^{\mathrm{f}}_{N-k,\theta_k\omega} = \sum_{t=0}^{N-k} \mathcal{Z}_{[k],[t+k]}(t) N^{1/2} P(N-k-t) \exp(\widetilde{\Phi}_{[t+k],[N]}(N-k-t)).$$

By (3.24), we then obtain

$$(3.25) \quad \exists \lim_{\substack{N\to\infty \\ N\in\eta}} N^{1/2} Z^{\mathrm{f}}_{N-k,\theta_k\omega} = \sum_{t=0}^{\infty} \mathcal{Z}_{[k],[t+k]}(t) \widetilde{L}_{[t+k],\eta} = [(1-B)^{-1}\widetilde{L}]_{[k],\eta}$$

since

$$(3.26) \quad \sum_{t=0}^{\infty} \mathcal{Z}_{\alpha,\gamma}(t) = \sum_{t=0}^{\infty}\sum_{k=0}^{\infty} M^{*k}_{\alpha,\gamma}(t) = \sum_{k=0}^{\infty} B^{*k}_{\alpha,\gamma} = [(I-B)^{-1}]_{\alpha,\gamma}.$$

*The critical case.* For $N \in \eta$ and $k \leq N$,

$$Z^{\mathrm{f}}_{N-k,\theta_k\omega} = \sum_{\gamma}\sum_{t=0}^{N-k} \mathcal{Z}_{[k],\gamma}(t) P(N-k-t) \exp(\widetilde{\Phi}_{\gamma,\eta}(N-k-t)).$$

By the previous results and using (3.24), we obtain that for every $k \in \mathbb{N}$,

$$
\begin{aligned}
\exists \lim_{N\to\infty, N\in\eta} Z^{\mathrm{f}}_{N-k,\theta_k\omega} &= \xi_{[k]} \frac{T}{2\pi} \frac{\langle \zeta, \widetilde{L}_{\cdot,\eta}\rangle}{\langle \zeta, L\xi\rangle} \int_0^1 \frac{dt}{t^{1/2}(1-t)^{1/2}} \\
&= \xi_{[k]} \frac{T}{2} \frac{\langle \zeta, \widetilde{L}_{\cdot,\eta}\rangle}{\langle \zeta, L\xi\rangle},
\end{aligned}
\tag{3.27}
$$

where we denote the canonical scalar product in $\mathbb{R}^{\mathbb{S}}$ by $\langle\cdot,\cdot\rangle$:

$$\langle \varphi, \psi \rangle := \sum_{\alpha\in\mathbb{S}} \varphi_\alpha \psi_\alpha, \qquad \varphi, \psi \in \mathbb{R}^{\mathbb{S}}.$$

*The localized case.* By (3.22),

$$e^{-F_\omega N} Z^f_{N-k,\theta_k\omega}$$

$$= e^{-F_\omega N} \sum_{t=0}^{N-k} \mathcal{Z}_{[k],[N-t]}(N-k-t) P(t) \exp(\widetilde{\Phi}_{[N-t],[N]}(t))$$

$$= e^{-F_\omega k} \sum_{\gamma \in \mathbb{S}} \sum_{t=0}^{N-k} e^{-F_\omega t} P(t) [\exp(\widetilde{\Phi}_{\gamma,[N]}(t)) e^{-F_\omega(N-k-t)} \mathcal{Z}_{[k],\gamma}(N-k-t)].$$

Since, by (3.12), the expression in brackets converges as $N \to \infty$ and $N \in [t] + \gamma$, we obtain

$$\exists \lim_{\substack{N \to \infty \\ N \in \eta}} e^{-F_\omega N} Z^f_{N-k,\theta_k\omega} = \xi_{[k]} e^{-F_\omega k} \left( \frac{T}{\mu} \sum_{\gamma \in \mathbb{S}} \sum_{t=0}^{\infty} e^{-F_\omega t} P(t) \exp(\widetilde{\Phi}_{\gamma,\eta}(t)) \zeta_\gamma \right).$$

Observe that the term in parentheses is just a function of $\eta$. $\square$

**4. A preliminary analysis of the polymer measure.** In this section, we give some preliminary material which will be used in Section 5 for the proof of the scaling limits of our models. We are going to show that the core of the polymer measure is encoded in its zero level set and that the law of the latter is expressible in terms of the partition function. This explains the crucial importance of the partition function for the study of $\mathbf{P}^a_{N,\omega}$.

We start by giving a very useful decomposition of $\mathbf{P}^a_{N,\omega}$. The intuitive idea is that a path $(S_n)_{n \leq N}$ can be split into two main ingredients:

- the family $(\tau_k)_{k=0,1,\ldots}$ of *returns to zero* of $S$ (defined in Section 2.2);
- the family of *excursions from zero* $(S_{i+\tau_{k-1}} : 0 \leq i \leq \tau_k - \tau_{k-1})_{k=1,2,\ldots}$.

Moreover, since each excursion can be either positive or negative, it is also useful to separately consider the signs of the excursions $\sigma_k := \text{sign}(S_{\tau_{k-1}+1})$ and the absolute values $(e_k(i) := |S_{i+\tau_{k-1}}| : i = 0, \ldots, \tau_k - \tau_{k-1})$. Observe that these are trivial for an excursion with length 1: in fact, if $\tau_k = \tau_{k-1} + 1$, then $\sigma_k = 0$ and $e_k(0) = e_k(1) = 0$.

Let us first consider the returns $(\tau_k)_{k \leq \iota_N}$ under $\mathbf{P}^a_{N,\omega}$, where we recall the definition $\iota_N = \sup\{k : \tau_k \leq N\}$. The law of this process can be viewed as a probability measure $p^a_{N,\omega}$ on the class of subsets of $\{1, \ldots, N\}$: indeed, for $A \subseteq \{1, \ldots, N\}$, writing

(4.1) $\quad A = \{t_1, \ldots, t_{|A|}\}, \qquad 0 =: t_0 < t_1 < \cdots < t_{|A|} \leq N,$

we can set

(4.2) $\quad\quad\quad p^a_{N,\omega}(A) := \mathbf{P}^a_{N,\omega}(\tau_i = t_i, i \leq \iota_N).$

From the definition (1.6) of $\mathbf{P}^a_{N,\omega}$ and from the strong Markov property of $\mathbf{P}$, we then have the following basic lemma:

LEMMA 4.1. *With the notation of (4.1), for $A \subseteq \{1, \ldots, N\}$ if $a = \mathrm{c}$, then $p^{\mathrm{c}}_{N,\omega}(A) \neq 0$ if and only if $t_{|A|} = N$ and, in this case,*

$$(4.3) \qquad p^{\mathrm{c}}_{N,\omega}(A) = \frac{1}{Z^{\mathrm{c}}_{N,\omega}} \prod_{i=1}^{|A|} M_{[t_{i-1}],[t_i]}(t_i - t_{i-1}),$$

*while for $a = \mathrm{f}$,*

$$(4.4) \qquad \begin{aligned} p^{\mathrm{f}}_{N,\omega}(A) &= \frac{1}{Z^{\mathrm{f}}_{N,\omega}} \left[ \prod_{i=1}^{|A|} M_{[t_{i-1}],[t_i]}(t_i - t_{i-1}) \right] \\ &\quad \times P(N - t_{|A|}) \exp(\widetilde{\Phi}_{[t_{|A|}],[N]}(N - t_{|A|})). \end{aligned}$$

Thus the law of the zero level set is explicitly given in terms of the kernel $M_{\alpha,\beta}(n)$ and the partition function $Z^a_{N,\omega}$. The following two lemmas [that again follow from definition (1.6)] show that, conditionally on the zero level set, the signs are independent and the excursions are just the excursions of the unperturbed random walk $S$ under $\mathbf{P}$. This shows that the zero level set is indeed the core of the polymer measure $\mathbf{P}^a_{N,\omega}$.

LEMMA 4.2. *Conditionally on $\{\iota_N, (\tau_j)_{j \leq \iota_N}\}$, under $\mathbf{P}^a_{N,\omega}$, the signs $(\sigma_k)_{k \leq \iota_N + 1}$ form an independent family. For $k \leq \iota_N$, the conditional law of $\sigma_k$ is specified as follows:*

- *if $\tau_k = 1 + \tau_{k-1}$, then $\sigma_k = 0$;*
- *if $\tau_k > 1 + \tau_{k-1}$, then $\sigma_k$ can take the two values $\pm 1$, with*

$$(4.5) \quad \mathbf{P}^a_{N,\omega}(\sigma_k = +1 | \iota_N, (\tau_j)_{j \leq \iota_N}) = \frac{1}{1 + \exp\{-(\tau_k - \tau_{k-1})h_\omega + \Sigma_{[\tau_{k-1}],[\tau_k]}\}}.$$

*For $a = \mathrm{f}$, when $\tau_{\iota_N} < N$, there is a last incomplete excursion in the interval $\{0, \ldots, N\}$, whose sign $\sigma_{\iota_N+1}$ is also specified by (4.5), provided we set $\tau_{\iota_N+1} := N$.*

LEMMA 4.3. *Conditionally on $\{\iota_N, (\tau_j)_{j \leq \iota_N}, (\sigma_j)_{j \leq \iota_N+1}\}$, the excursions $(e_k(\cdot))_{k \leq \iota_N+1}$ form an independent family under $\mathbf{P}^a_{N,\omega}$. For $k \leq \iota_N$, the conditional law of $e_k(\cdot)$ on the event $\{\tau_{k-1} = \ell_0, \tau_k = \ell_1\}$ is specified for $f = (f_i)_{i=0,\ldots,\ell_1-\ell_0}$ by*

$$(4.6) \quad \begin{aligned} \mathbf{P}^a_{N,\omega}&(e_k(\cdot) = f | \iota_N, (\tau_j)_{j \leq \iota_N}, (\sigma_j)_{j \leq \iota_N+1}) \\ &= \mathbf{P}(S_i = f_i : i = 0, \ldots, \ell_1 - \ell_0 | S_i > 0 : i = 1, \ldots, \ell_1 - \ell_0 - 1, \\ &\qquad S_{\ell_1 - \ell_0} = 0). \end{aligned}$$

For $a = \mathrm{f}$, when $\tau_{\iota_N} < N$, the conditional law on the event $\{\tau_{\iota_N} = \ell < N\}$ of the last incomplete excursion $e_{\iota_N+1}(\cdot)$ is specified for $f = (f_i)_{i=0,\ldots,N-\ell}$ by

$$
\begin{aligned}
\mathbf{P}^a_{N,\omega}(e_{\iota_N+1}(\cdot) &= f | \iota_N, (\tau_j)_{j \le \iota_N}, (\sigma_j)_{j \le \iota_N+1}) \\
&= \mathbf{P}(S_i = f_i : i = 0, \ldots, N - \ell | S_i > 0 : i = 1, \ldots, N - \ell).
\end{aligned}
\tag{4.7}
$$

We stress that Lemmas 4.1, 4.2 and 4.3 fully characterize the polymer measure $\mathbf{P}^a_{N,\omega}$. It is worth stressing that, conditionally on $(\tau_k)_{k \in \mathbb{N}}$, the joint distribution of $(\sigma_j, e_j)_{j \le \iota_N}$ *does not depend on $N$*. In this sense, all the $N$-dependence is contained in the law $p^a_{N,\omega}$ of the zero level set. This fact will be exploited in the next section.

**5. Proof of Theorem 1.3 and Proposition 2.1.** In this section, we show that the measures $\mathbf{P}^a_{N,\omega}$ converge under Brownian rescaling, proving Theorem 1.3 and Proposition 2.1. The results and proofs closely follow those of [9] and we shall refer to this paper for several technical lemmas; for the tightness of $(Q^a_{N,\omega})_{N \in \mathbb{N}}$ in $C([0,1])$, we refer to [8].

LEMMA 5.1. *For any $\omega$ and $a = \mathrm{c}, \mathrm{f}$, the sequence $(Q^a_{N,\omega})_{N \in \mathbb{N}}$ is tight in $C([0,1])$.*

Hereafter, we separately consider the three regimes $\delta^\omega > 1$, $\delta^\omega < 1$ and $\delta^\omega = 1$.

5.1. *The localized regime ($\delta^\omega > 1$).* We prove point (3) of Theorem 1.3. By Lemma 5.1, it is enough to prove that $\mathbf{P}^a_{N,\omega}(|X^N_t| > \varepsilon) \to 0$ for all $\varepsilon > 0$ and $t \in [0,1]$ and one can obtain this estimate explicitly. We point out, however, that in this regime, one can avoid using the compactness lemma and can obtain a stronger result by elementary means: observe that for any $k, n \in \mathbb{N}$ such that $n > 1$ and $k + n \le N$, we have

$$
\begin{aligned}
\mathbf{P}^a_{N,\omega}&(S_k = S_{k+n} = 0, S_{k+i} \ne 0 \text{ for } i = 1, \ldots, n-1) \\
&\le \frac{1/2(1 + \exp(\sum_{i=1}^n (\omega^{(-1)}_{k+i} - \omega^{(+1)}_{k+i})))}{Z^{\mathrm{c}}_{n,\theta_k \omega}} =: \widehat{K}_k(n),
\end{aligned}
\tag{5.1}
$$

and that this holds both for $a = \mathrm{c}$ and $a = \mathrm{f}$. Inequality (5.1) is obtained by using the Markov property of $S$ both in the numerator and the denominator of the expression (1.6) defining $\mathbf{P}^a_{N,\omega}(\cdot)$ after having bounded $Z^a_{N,\omega}$ from below by inserting the event $S_k = S_{k+n} = 0$. Of course, $\lim_{n \to \infty}(1/n) \log \widehat{K}_k(n) = -\mathrm{F}_\omega$ uniformly in $k$ [notice that $\widehat{K}_{k+T}(n) = \widehat{K}_k(n)$].

Therefore, if we fix $\varepsilon > 0$, by the union bound, we obtain (recalling that $\{\tau_j\}_j$ and $\iota_N$ were defined in Section 2.2), for some $c > 0$,

$$\mathbf{P}^a_{N,\omega}\left(\max_{j=1,2,\ldots,\iota_N} \tau_j - \tau_{j-1} > (1+\varepsilon)\log N/\mathrm{F}_\omega\right)$$

$$\leq \sum_{k \leq N-(1+\varepsilon)\log N/\mathrm{F}_\omega} \sum_{n > (1+\varepsilon)\log N/\mathrm{F}_\omega} \widehat{K}_k(n)$$

$$\leq N \sum_{n > (1+\varepsilon)\log N/\mathrm{F}_\omega} \max_{k=0,\ldots,T-1} \widehat{K}_k(n)$$

$$\leq \frac{c}{N^\varepsilon}.$$

Let us start with the constrained case. Notice that $\mathbf{P}^{\mathrm{c}}_{N,\omega}(\mathrm{d}S)$-a.s. we have $\tau_{\iota_N} = N$ and hence $\max_{j \leq \iota_N} \tau_j - \tau_{j-1} \geq \max_{n=1,\ldots,N} |S_n|$, since $|S_{n+1} - S_n| \leq 1$. We then immediately obtain that for any $C > 1/\mathrm{F}_\omega$,

(5.2) $$\lim_{N \to \infty} \mathbf{P}^{\mathrm{c}}_{N,\omega}\left(\max_{n=1,\ldots,N} |S_n| > C \log N\right) = 0,$$

which is, of course, a much stronger statement than the scaling limit of point (3) of Theorem 1.3. If we instead consider the measure $\mathbf{P}^{\mathrm{f}}_{N,\omega}$, the length of the last excursion must also be taken into account; however, an argument very close to the one used in (5.1) also yields that the last excursion is exponentially bounded (with the same exponent) and the proof of point (3) of Theorem 1.3 is complete.

5.2. *The strictly delocalized regime* ($\delta^\omega < 1$). We prove point (1) of Theorem 1.3 and Proposition 2.1. We set, for $t \in \{1,\ldots,N\}$,

$$D_t := \inf\{k = 1,\ldots,N : k > t, S_k = 0\},$$
$$G_t := \sup\{k = 1,\ldots,N : k \leq t, S_k = 0\}.$$

The following result shows that in the strictly delocalized regime, as $N \to \infty$, the visits to zero under $\mathbf{P}^a_{N,\omega}$ tend to be very few and concentrated at a finite distance from the origin if $a = \mathrm{f}$ and from $0$ or $N$ if $a = \mathrm{c}$.

LEMMA 5.2. *If $\delta^\omega < 1$, there exists a constant $C > 0$ such that for all $L > 0$,*

$$\limsup_{N \to \infty}[\mathbf{P}^{\mathrm{f}}_{N,\omega}(G_N \geq L) + \mathbf{P}^{\mathrm{c}}_{N,\omega}(G_{N/2} \geq L) + \mathbf{P}^{\mathrm{c}}_{N,\omega}(D_{N/2} \leq N - L)]$$

$$\leq CL^{-1/2}.$$

PROOF. We consider, for example, $\mathbf{P}^{\mathrm{c}}_{N,\omega}(G_{N/2} \geq L)$. Using Lemma 4.1 to compute this probability, recalling definition (3.5) of the kernel $\mathcal{Z}_{\alpha,\beta}(n)$ and using (3.18), we obtain

$$\mathbf{P}^{\mathrm{c}}_{N,\omega}(G_{N/2} \geq L) = \sum_{x=L}^{\lfloor N/2 \rfloor} \mathcal{Z}_{0,[x]}(x) \sum_{z=\lfloor N/2 \rfloor+1}^{N} \frac{M_{[x],[z]}(z-x)\mathcal{Z}_{[z],[N]}(N-z)}{\mathcal{Z}_{0,[N]}(N)}$$

$$\leq C_1 N^{3/2} \sum_{x=L}^{\lfloor N/2 \rfloor} x^{-3/2} \sum_{z=\lfloor N/2 \rfloor+1}^{N} (z-x)^{-3/2}(N+1-z)^{-3/2}$$

$$\leq C_2 L^{-1/2}$$

for some positive constants $C_1$ and $C_2$, and the proof is complete. □

*The signs.* Since the zeros are concentrated near the boundary, to complete the proof it is enough to argue as in the proof of Theorem 9 in [9]. More precisely, by Lemma 5.2, for large $N$, the typical paths of $\mathbf{P}^{a}_{N,\omega}$ are essentially made up of one big excursion whose absolute value converges in law to the Brownian excursion $\{e_t\}_{t \in [0,1]}$ for $a = \mathrm{c}$ and to the Brownian meander $\{m_t\}_{t \in [0,1]}$ for $a = \mathrm{f}$ by standard invariance principles (cf. [18] and [3]). Therefore, to complete the proof, we only have to show the existence the limit (as $N \to \infty$ along $[N] = \eta$) of the probability that the process (away from $\{0,1\}$) lives in the upper half-plane. In the general case, we have different limits, depending on the sequence $[N] = \eta$ and on $a = \mathrm{f}, \mathrm{c}$, while if $\omega \notin \mathcal{P}$, all such limits coincide.

We start with the constrained case. Given Lemma 5.2, it is sufficient to show that

(5.3) $$\exists \lim_{\substack{N \to \infty \\ N \in \eta}} \mathbf{P}^{\mathrm{c}}_{N,\omega}(S_{N/2} > 0) =: \mathrm{p}^{\mathrm{c}}_{\omega,\eta}.$$

Formula (5.3) follows from the fact that

$$\mathbf{P}^{\mathrm{c}}_{N,\omega}(S_{N/2} > 0)$$

$$= \sum_{\alpha,\beta \in \mathbb{S}} \sum_{x<N/2} \sum_{y>N/2} \frac{\mathcal{Z}_{0,\alpha}(x)\rho^{+}_{\alpha,\beta}(y-x)M_{\alpha,\beta}(y-x)\mathcal{Z}_{\beta,[N]}(N-y)}{\mathcal{Z}_{0,[N]}(N)},$$

where for all $z \in \mathbb{N}$ and $\alpha, \beta \in \mathbb{S}$, we set

(5.4) $$\rho^{+}_{\alpha,\beta}(z) := \frac{1}{1+\exp(-zh_\omega + \Sigma_{\alpha,\beta})};$$

see (4.5). By the dominated convergence theorem and by (3.10) and (3.26), we get

$$\exists \lim_{\substack{N \to \infty \\ N \in \eta}} N^{3/2} \sum_{x<N/2} \sum_{y>N/2} \mathcal{Z}_{0,\alpha}(x)\rho^{+}_{\alpha,\beta}(y-x)M_{\alpha,\beta}(y-x)\mathcal{Z}_{\beta,\eta}(N-y)$$

$$= [(1-B)^{-1}]_{0,\alpha} c_K \tfrac{1}{2} \exp(\omega_\beta^{(0)}) [(1-B)^{-1}]_{\beta,\eta}.$$

By (3.18), it then follows that (5.3) holds true with

$$(5.5) \quad \mathtt{p}_{\omega,\eta}^{\mathrm{c}} := \frac{\sum_{\alpha,\beta}[(1-B)^{-1}]_{0,\alpha} c_K (1/2) \exp(\omega_\beta^{(0)}) [(1-B)^{-1}]_{\beta,\eta}}{[(1-B)^{-1} L (1-B)^{-1}]_{0,\eta}}.$$

REMARK 5.3. Observe that by (3.10):

- if $h_\omega > 0$, then in (5.5), the denominator is equal to the numerator, so that $\mathtt{p}_{\omega,\eta}^{\mathrm{c}} = 1$ for all $\eta$;
- if $h_\omega = 0$ and $\Sigma \equiv 0$, then in (5.5), the denominator is equal to *twice* the numerator, so that $\mathtt{p}_{\omega,\eta}^{\mathrm{c}} = 1/2$ for all $\eta$;
- in the remaining case, that is, if $\omega \in \mathcal{P}$, then in general, $\mathtt{p}_{\omega,\eta}^{\mathrm{c}}$ depends on $\eta$.

Let us now consider the free case. This time, it is sufficient to show that

$$(5.6) \quad \exists \lim_{\substack{N \to \infty \\ N \in \eta}} \mathbf{P}_{N,\omega}^{\mathrm{f}}(S_N > 0) =: \mathtt{p}_{\omega,\eta}^{\mathrm{f}}.$$

However, we can write

$$\mathbf{P}_{N,\omega}^{\mathrm{f}}(S_N > 0) = \sum_\alpha \sum_{x<N} \frac{\mathcal{Z}_{0,\alpha}(x) \cdot (1/2) P(N-k)}{Z_{N,\omega}^{\mathrm{f}}}$$

and by using (3.22), (3.26) and (3.24), we obtain that (5.6) holds with

$$(5.7) \quad \mathtt{p}_{\omega,\eta}^{\mathrm{f}} = \frac{\sum_\alpha [(1-B)^{-1}]_{0,\alpha} c_K}{[(1-B)^{-1} \widetilde{L}]_{0,\eta}}.$$

REMARK 5.4. Again, observe that by (3.24):

- if $h_\omega > 0$, then in (5.7), the denominator is equal to the numerator and $\mathtt{p}_{\omega,\eta}^{\mathrm{f}} = 1$ for all $\eta$;
- if $h_\omega = 0$ and $\Sigma \equiv 0$, then in (5.7), the denominator is equal to *twice* the numerator, so that $\mathtt{p}_{\omega,\eta}^{\mathrm{f}} = 1/2$ for all $\eta$;
- in the remaining case, that is, if $\omega \in \mathcal{P}$, the in general, $\mathtt{p}_{\omega,\eta}^{\mathrm{f}}$ depends on $\eta$ and is different from $\mathtt{p}_{\omega,\eta}^{\mathrm{c}}$.

The proof of point (1) of Theorem 1.3 and Proposition 2.1 is then concluded.

5.3. *The critical regime* ($\delta^\omega = 1$). We prove point (2) of Theorem 1.3 and Proposition 2.1. As in the previous section, we first determine the asymptotic behavior of the zero level set and then pass to the study of the signs of the excursions.

*The zero level set.* We introduce the random closed subset $\mathcal{A}_N^a$ of $[0,1]$ describing the zero set of the polymer of size $N$ rescaled by a factor $1/N$:

$$\mathbb{P}(\mathcal{A}_N^a = A/N) = p_{N,\omega}^a(A), \qquad A \subseteq \{0,\ldots,N\};$$

recall (4.2). Let us denote by $\mathcal{F}$ the class of *all closed subsets* of $\mathbb{R}^+ := [0,+\infty)$. We are going to impose on $\mathcal{F}$ a topological and measurable structure so that we can view the law of $\mathcal{A}_N^a$ as a probability measure on $\mathcal{F}$ and can study the weak convergence of $\mathcal{A}_N^a$.

We endow $\mathcal{F}$ with the topology of Matheron (cf. [21] and [12], Section 3) which is a metrizable topology. To define it, we associate to a closed subset $F \subseteq \mathbb{R}^+$ the *compact nonempty* subset $\widetilde{F}$ of the interval $[0,\pi/2]$ defined by $\widetilde{F} := \arctan(F \cup \{+\infty\})$. The metric $\rho(\cdot,\cdot)$ we introduce on $\mathcal{F}$ is then

$$(5.8) \qquad \rho(F, F') := \max\left\{\sup_{t \in \widetilde{F}} d(t, \widetilde{F'}), \sup_{t' \in \widetilde{F'}} d(t', \widetilde{F})\right\}, \qquad F, F' \in \mathcal{F},$$

where $d(s, A) := \inf\{|t-s|, t \in A\}$ is the standard distance between a point and a set. The r.h.s. of (5.8) is the so-called Hausdorff metric between the compact nonempty sets $\widetilde{F}, \widetilde{F'}$.

Thus, by definition, a sequence $\{F_n\}_n \subset \mathcal{F}$ converges to $F \in \mathcal{F}$ if and only if $\rho(F_n, F) \to 0$. This is equivalent to requiring that for each open set $G$ and each compact $K \subset \mathbb{R}^+$,

$$(5.9) \qquad \begin{aligned} F \cap G \neq \varnothing &\implies F_n \cap G \neq \varnothing \text{ eventually,} \\ F \cap K = \varnothing &\implies F_n \cap K = \varnothing \text{ eventually.} \end{aligned}$$

Another necessary and sufficient condition for $F_n \to F$ is that $d(t, F_n) \to d(t, F)$ for every $t \in \mathbb{R}^+$.

This topology makes $\mathcal{F}$ a separable and compact metric space [21], Theorem 1-2-1, in particular, a Polish space. Endowing $\mathcal{F}$ with the Borel $\sigma$-field, we have that the space $\mathcal{M}_1(\mathcal{F})$ of probability measures on $\mathcal{F}$ is compact with the topology of weak convergence.

The crucial result is the convergence in distribution as $N \to \infty$ of the random set $\mathcal{A}_N^a$ toward the zero set of a Brownian motion for $a = \mathrm{f}$ or of a Brownian bridge for $a = \mathrm{c}$.

PROPOSITION 5.5. *If $\delta^\omega = 1$, then as $N \to \infty$,*

$$(5.10) \qquad \mathcal{A}_N^{\mathrm{f}} \implies \{t \in [0,1] : B(t) = 0\},$$

$$(5.11) \qquad \mathcal{A}_N^{\mathrm{c}} \implies \{t \in [0,1] : \beta(t) = 0\}.$$

The proof of Proposition 5.5 is achieved by comparing the law of $\mathcal{A}_N^{\mathrm{f}}$ and $\mathcal{A}_N^{\mathrm{c}}$ with the law of a random set $\mathcal{R}_N$ defined as follows. With the notation introduced in Section 3.1, we introduce the rescaled random set $\mathcal{R}_N$:

$$\mathcal{R}_N := \mathrm{range}\{\tau_i/N, i \geq 0\} = \tau/N \subset \mathbb{R}^+$$

under $\mathbb{P}_{[0]}$. Notice that for any $A = \{t_1, \ldots, t_{|A|}\} \subset \{1, \ldots, N\}$, we have (setting $t_0 := 0$)

$$
\begin{aligned}
&\mathbb{P}_{[0]}(\tau \cap \{1, \ldots, N\} = A) \\
&\qquad = \left[\prod_{i=1}^{|A|} M_{[t_{i-1}],[t_i]}(t_i - t_{i-1})\right] Q_{t_{|A|}}(N - t_{|A|}) \frac{\xi_{[t_{|A|}]}}{\xi_{[0]}},
\end{aligned}
\tag{5.12}
$$

where $Q_\alpha(t) := \sum_\beta \sum_{s=t+1}^\infty \Gamma_{\alpha,\beta}(s)$.

The key step in proving Proposition 5.5 is given by the following lemma whose proof uses the theory of *regenerative sets* and their connection with *subordinators*; see [12].

LEMMA 5.6. *The random set $\mathcal{R}_N$ converges in distribution to $\{t \geq 0 : B(t) = 0\}$.*

PROOF. Recalling the definition (3.8) of $\tau^\beta$, we introduce the random set

$$\mathcal{R}_N^\beta := \operatorname{range}\{\tau_k/N : k \geq 0, J_k = \beta\} = \tau^\beta/N, \qquad \beta \in \mathbb{S},$$

under $\mathbb{P}_{[0]}$. Notice that $\mathcal{R}_N = \bigcup_\beta \mathcal{R}_N^\beta$. We divide the rest of the proof into two steps.

*Step* 1. This is the main step: we prove that the law of $\mathcal{R}_N^\beta$ converges to the law of $\{t \geq 0 : B(t) = 0\}$. For this, we follow the proof of Lemma 5 in [9]. Let $\{P(t)\}_{t\geq 0}$ be a Poisson process with rate $\gamma > 0$, independent of $(T_i^\beta)_{i\geq 0}$. Then $\sigma_t = [T_1^\beta + \cdots + T_{P(t)}^\beta]/N$ is a nondecreasing right-continuous process with independent stationary increments and $\sigma_0 = 0$, that is, $\sigma = (\sigma_t)_{t\geq 0}$ is a subordinator. By the standard theory of Lévy processes, the law of $\sigma$ is characterized by the Laplace transform of its one-time distributions,

$$\mathbb{E}[\exp(-\lambda \sigma_t)] = \exp(-t\phi_N(\lambda)), \qquad \lambda \geq 0, t \geq 0,$$

for a suitable function $\phi_N : [0, \infty) \mapsto [0, \infty)$, called the Lévy exponent, which has a canonical representation, the Lévy–Khinchine formula (see, e.g., (1.15) in [12]):

$$\phi_N(\lambda) = \int_{(0,\infty)} (1 - e^{-\lambda s})\gamma \mathbb{P}(T_1^\beta/N \in ds) = \gamma \sum_{n=1}^\infty (1 - \exp(-\lambda n/N))q^\beta(n).$$

We denote the closed range $\{\sigma_t : t \geq 0\}$ of the subordinator $\sigma$ by $\widehat{\mathcal{R}}_N^\beta$. Then, following [12], $\widehat{\mathcal{R}}_N^\beta$ is a regenerative set. Moreover, $\mathcal{R}_N^\beta = T_0^\beta/N + \widehat{\mathcal{R}}_N^\beta$.

Notice now that the law of the regenerative set $\widehat{\mathcal{R}}_N^\beta$ is invariant under the change of time scale $\sigma_t \longrightarrow \sigma_{ct}$ for $c > 0$ and, in particular, independent of $\gamma > 0$. Since $\phi_N \longrightarrow c\phi_N$ under this change of scale, we can fix $\gamma = \gamma_N$ such that $\phi_N(1) = 1$ and this will hereafter be implicitly assumed. Then, by Proposition (1.14) of [12], the law of $\widehat{\mathcal{R}}_N^\beta$ is uniquely determined by $\phi_N$.

The asymptotic behavior of $q^\beta$ given in (3.15) easily yields $\phi_N(\lambda) \to \lambda^{1/2} =: \Phi_{BM}(\lambda)$ as $N \to \infty$. It is now a matter of applying the result in [12] Section 3, to obtain that $\widehat{\mathcal{R}}_N^\beta$ converges in law to the regenerative set corresponding to $\Phi_{BM}$. However, by direct computation, one obtains that the latter is nothing but the zero level set of a Brownian motion, therefore $\widehat{\mathcal{R}}_N^\beta \Rightarrow \{t \in [0,1] : B(t) = 0\}$. From the fact that $T_0^\beta/N \to 0$ a.s., the same weak convergence for $\mathcal{R}_N^\beta$ follows immediately.

*Step* 2. Notice that $\mathcal{R}_N = \bigcup_\beta \mathcal{R}_N^\beta$ is the union of nonindependent sets. Therefore, although we know that each $\mathcal{R}_N^\beta$ converges in law to $\{t \geq 0 : B(t) = 0\}$, it is not trivial that $\mathcal{R}_N$ converges to the same limit. We start by showing that for every positive $t \geq 0$, the distance between the first point in $\mathcal{R}_N^\alpha$ after $t$ and the first point in $\mathcal{R}_N^\beta$ after $t$ converges to zero in probability. More precisely, for any closed set $F \subset [0,\infty)$, we set

$$(5.13) \qquad d_t(F) := \inf(F \cap (t,\infty))$$

and we claim that for all $\alpha, \beta \in \mathbb{S}$ and $t \geq 0$, $|d_t(\mathcal{R}_N^\alpha) - d_t(\mathcal{R}_N^\beta)| \to 0$ in probability.

Recalling (3.14) and the notation introduced there, we can write, for all $\epsilon > 0$,

$$\mathbb{P}_{[0]}(d_t(\mathcal{R}_N^\alpha) \geq d_t(\mathcal{R}_N^\beta) + \epsilon)$$
$$= \sum_\gamma \sum_{y=0}^{\lfloor Nt \rfloor} U_{0,\gamma}(y) \sum_{z=\lfloor Nt \rfloor - y + 1}^\infty q^{(\gamma;\beta)}(z) \sum_{w=\lfloor N\epsilon \rfloor}^\infty q^{(\beta;\alpha)}(w).$$

Arguing as in the proof of (3.15), it is easy to obtain the bound $q^{(\beta;\alpha)}(w) \leq C_1 w^{-3/2}$ and by (3.16), we have $U_{0,\gamma}(y) \leq C_2 y^{-1/2}$, with $C_1, C_2$ positive constants. Therefore,

$$\mathbb{P}_{[0]}(d_t(\mathcal{R}_N^\alpha) \geq d_t(\mathcal{R}_N^\beta) + \epsilon)$$
$$\leq \frac{C_3}{N^{1/2}} \left( \int_0^{t/T} dy \int_{(t-y)/T}^\infty dz \int_{\epsilon/T}^\infty dw \frac{1}{y^{1/2} z^{3/2} w^{3/2}} \right)$$

for some positive constant $C_3$, having used the convergence of the Riemann sums to the corresponding integral. The same computations can be performed where $\alpha$ is exchanged with $\beta$, hence the claim is proved.

Now, notice that $d_t(\mathcal{R}_N) = \min_{\alpha \in \mathbb{S}} d_t(\mathcal{R}_N^\alpha)$ and since $\mathbb{S}$ is a finite set, we also have that $|d_t(\mathcal{R}_N) - d_t(\mathcal{R}_N^\beta)| \to 0$ in probability for any fixed $\beta \in \mathbb{S}$. Since we already know that $\mathcal{R}_N^\beta$ converges weakly to the law of $\{t \geq 0 : B(t) =$ 0}, the analogous statement for $\mathcal{R}_N$ follows by standard arguments. More precisely, let us look at $(\mathcal{R}_N, \mathcal{R}_N^\beta)$ as a random element of the space $\mathcal{F} \times \mathcal{F}$: by the compactness of $\mathcal{F}$, it suffices to take any convergent subsequence $(\mathcal{R}_{k_n}, \mathcal{R}_{k_n}^\beta) \Rightarrow (\mathfrak{B}, \mathfrak{C})$ and show that $\mathbb{P}(\mathfrak{B} \neq \mathfrak{C}) = 0$. However, we can write

$$\{\mathfrak{B} \neq \mathfrak{C}\} = \bigcup_{t \in \mathbb{Q}^+} \bigcup_{n \in \mathbb{N}} \{|d_t(\mathfrak{B}) - d_t(\mathfrak{C})| > 1/n\}$$

by the right-continuity of $t \mapsto d_t$ and, by the portmanteau theorem, we have

$$\mathbb{P}_{[0]}(|d_t(\mathfrak{B}) - d_t(\mathfrak{C})| > 1/n) \leq \limsup_{N \to \infty} \mathbb{P}_{[0]}(|d_t(\mathcal{R}_N) - d_t(\mathcal{R}_N^\beta)| > 1/n) = 0$$

because $|d_t(\mathcal{R}_N) - d_t(\mathcal{R}_N^\beta)| \to 0$ in probability. $\square$

PROOF OF PROPOSITION 5.5: EQUATION (5.11). First, we compute the Radon–Nikodym derivative $f_N^c$ of the law of $\mathcal{A}_N^c \cap [0, 1/2]$ with respect to the law of $\mathcal{R}_N^{1/2} := \mathcal{R}_N \cap [0, 1/2]$, using (4.3) and (5.12). For $F = \{t_1/N, \ldots, t_k/N\} \subset [0, 1/2]$ with $0 =: t_0 < t_1 < \cdots < t_k$ integer values, the value of $f_N^c$ at $\mathcal{R}_N^{1/2} = F$ depends only on $g_{1/2}(F)$ and is given by

$$f_N^c(g_{1/2}(F)) = f_N^c(t_k/N) = \frac{\sum_{n=N/2}^{N} M_{[t_k],[n]}(n - t_k) \mathcal{Z}_{[n],[N]}(N - n)}{\mathcal{Z}_{0,[N]}(N) Q_{[t_k]}(N/2 - t_k)} \frac{\xi_0}{\xi_{[t_k]}},$$

where for any closed set $F \subset [0, \infty)$, we set

(5.14) $$g_t(F) := \sup(F \cap [0, t]).$$

By (3.17), for all $\varepsilon > 0$ and uniformly in $g \in [0, 1/2 - \varepsilon]$,

$$f_N^c(g) \sim \frac{[L\xi]_{[Ng]}(T^2/(2\pi))(\zeta_{[N]}/\langle \zeta, L\xi \rangle)T^{-1} \int_0^{1/2} y^{-1/2}(1 - y - g)^{-3/2}\, dy}{(T^2/(2\pi))(\xi_0 \zeta_{[N]} \langle \zeta, L\xi \rangle)T^{-1}([L\xi]_{[Ng]}/\xi_{[Ng]})2(1/2 - g)^{-1/2}}$$

$$\times \frac{\xi_0}{\xi_{[Ng]}}$$

$$= \frac{\sqrt{1/2}}{1 - g} =: r(g).$$

If $\Psi$ is a bounded continuous functional on $\mathcal{F}$ such that $\Psi(F) = \Psi(F \cap [0, 1/2])$ for all $F \in \mathcal{F}$, then setting $Z_B := \{t \in [0,1] : B(t) = 0\}$ and $Z_\beta := \{t \in [0,1] : \beta(t) = 0\}$, we get

$$\mathbb{E}[\Psi(Z_\beta)] = \mathbb{E}[\Psi(Z_B) r(g_{1/2}(Z_B))];$$

see formula (49) in [9]. By Lemma 5.6 and the asymptotic behavior of $f_N^c$, we obtain

$$\mathbb{E}[\Psi(\mathcal{A}_N^c)] = \mathbb{E}[\Psi(\mathcal{R}_N^{1/2}) f_N^c(g_{1/2}(\mathcal{R}_N^{1/2}))] \stackrel{N \to \infty}{\longrightarrow} \mathbb{E}[\Psi(Z_B) r(g_{1/2}(Z_B))]$$

$$= \mathbb{E}[\Psi(Z_\beta)],$$

that is, $\mathcal{A}_N^c \cap [0, 1/2]$ converges to $Z_\beta \cap [0, 1/2]$. Notice now that the distribution of the random set $\{1 - t : t \in \mathcal{A}_N^c \cap [1/2, 1]\}$ under $\mathbf{P}_{N,\omega}^c$ is the same as the distribution of $\mathcal{A}_N^c \cap [0, 1/2]$ under $\mathbf{P}_{N,\overline{\omega}}^c$, where $\overline{\omega}_{[i]} := \omega_{[N-i]}$. Therefore, we obtain that $\mathcal{A}_N^c \cap [1/2, 1]$ converges to $Z_\beta \cap [0, 1/2]$ and the proof is complete.

PROOF OF PROPOSITION 5.5: EQUATION (5.10). By conditioning on the last zero, from (4.3) and (4.4), we see that if $\Psi$ is a bounded continuous functional on $\mathcal{F}$, then

$$\mathbb{E}[\Psi(\mathcal{A}_N^f)] = \sum_{k=0}^N \mathbb{E}\left[\Psi\left(\frac{k}{N}\mathcal{A}_k^c\right)\right] \frac{Z_{k,\omega}^c}{Z_{N,\omega}^f} P(N-k) \exp(\widetilde{\Phi}_{[k],[N]}(N-k)).$$

We denote by $\beta^t$ a Brownian bridge over the interval $[0, t]$, that is, a Brownian motion over $[0, t]$ conditioned to be 0 at time $t$, and we set $Z_{\beta^t} := \{s \in [0, t] : \beta^t(s) = 0\} \stackrel{d}{=} tZ_\beta$. By (5.11), it follows that if $k/N \to t$, then the random set $\frac{k}{N}\mathcal{A}_k^c$ converges in distribution to $Z_{\beta^t}$. Then, applying (3.17) and (3.27), we obtain, as $N \to \infty$ along $[N] = \eta$,

$$\begin{aligned}
\mathbb{E}[\Psi(\mathcal{A}_k^c)] &= \sum_{k=0}^N \sum_\gamma 1_{(k \in \gamma)} \mathbb{E}\left[\Psi\left(\frac{k}{N}\mathcal{A}_k^c\right)\right] \frac{Z_{k,\omega}^c}{Z_{N,\omega}^f} P(N-k) \exp(\widetilde{\Phi}_{\gamma,\eta}(N-k)) \\
&\longrightarrow \int_0^1 \mathbb{E}[\Psi(Z_{\beta^t})] \frac{1}{\pi t^{1/2}(1-t)^{1/2}} \, dt \\
&\quad \times \sum_\gamma \frac{1}{T} \frac{T^2}{2} \frac{\xi_0 \zeta_\gamma}{\langle \zeta, L\xi \rangle} \frac{\widetilde{L}_{\gamma,\eta}}{\xi_0(T/2)\langle \zeta, \widetilde{L}_{\cdot,\eta}\rangle / \langle \zeta, L\xi \rangle} \\
&= \mathbb{E}[\Psi(Z_B)].
\end{aligned}$$

Since the result does not depend on the subsequence $[N] = \eta$, we have indeed proven that $\mathcal{A}_N^f$ converges in distribution to $Z_B$. $\square$

*The signs.* In order to conclude the proof of point (2) of Theorem 1.3 and Proposition 2.1 in the critical case ($\delta^\omega = 1$), we closely follow the proof given in Section 8 of [9]. Having already proven the convergence of the zero level set, we only have to *paste* the excursions (recall Lemmas 4.2 and 4.3). The weak convergence under diffusive rescaling of $e_k(\cdot)$ for $k \leq \iota_N$ toward the Brownian excursion $e(\cdot)$ and of the last excursion $e_{\iota_N+1}(\cdot)$ for $a = f$ toward the Brownian meander $m(\cdot)$ has been proven in [18] and [3], respectively. It then only remains to focus on the signs.

We start with the constrained case: we are going to show that for all $t \in (0, 1)$,

(5.15) $$\exists \lim_{N \to \infty} \mathbf{P}_{N,\omega}^c(S_{\lfloor tN \rfloor} > 0) =: \mathtt{p}_\omega,$$

and that the limit is independent of $t$. We point out that we should actually fix the extremities of the excursion embracing $t$, that is, we should rather prove that

$$\lim_{N \to \infty} \mathbf{P}^{\mathrm{c}}_{N,\omega}(S_{\lfloor tN \rfloor} > 0 | G_{\lfloor tN \rfloor}/N \in (a-\varepsilon, a), D_{\lfloor tN \rfloor}/N \in (b, b+\varepsilon))$$
(5.16)
$$= \mathrm{p}_\omega,$$

for $a < t < b$ and $\varepsilon > 0$ (recall the definitions of $G_t$ and $D_t$ in Section 5.2), but in order to lighten the exposition, we will stick to (5.15), since proving (5.16) requires only minor changes.

We have, recalling (5.4),

$$\mathbf{P}^{\mathrm{c}}_{N,\omega}(S_{\lfloor tN \rfloor} > 0)$$

$$= \sum_{\alpha,\beta} \sum_{x < \lfloor tN \rfloor} \sum_{y > \lfloor tN \rfloor} \frac{\mathcal{Z}_{0,\alpha}(x) \rho^+_{\alpha,\beta}(y-x) M_{\alpha,\beta}(y-x) \mathcal{Z}_{\beta,[N]}(N-y)}{\mathcal{Z}_{0,[N]}(N)}.$$

By the dominated convergence theorem and by (3.17),

$$\exists \lim_{\substack{N \to \infty \\ N \in \eta}} N^{1/2} \sum_{x < \lfloor tN \rfloor} \sum_{y > \lfloor tN \rfloor} \mathcal{Z}_{0,\alpha}(x) \rho^+_{\alpha,\beta}(y-x) M_{\alpha,\beta}(y-x) \mathcal{Z}_{\beta,\eta}(N-y)$$

$$= \frac{1}{T^2} \int_0^t dx \int_t^1 dy \, [x(y-x)^3(1-y)]^{1/2} \left(\frac{T^2}{2\pi}\right)^2 \frac{\xi_0 \zeta_\alpha \xi_\beta \zeta_\eta}{\langle \zeta, L\xi \rangle^2} c_K \frac{1}{2} \exp(\omega^{(0)}_\beta);$$

see (3.10). We obtain that (5.15) holds with

(5.17)
$$\mathrm{p}_\omega := \frac{\sum_{\alpha,\beta} \zeta_\alpha c_K(1/2) \exp(\omega^{(0)}_\beta) \xi_\beta}{\langle \zeta, L\xi \rangle}.$$

REMARK 5.7. Observe the following. By (3.10):

- if $h_\omega > 0$, then in (5.17), the denominator is equal to the numerator so that $\mathrm{p}_\omega = 1$;
- if $h_\omega = 0$, and $\Sigma \equiv 0$ then in (5.17), the denominator is equal to *twice* the numerator so that $\mathrm{p}_\omega = 1/2$.

Let us now consider the free case. We are going to show that for all $t \in (0,1]$,

$$\lim_{\substack{N \to \infty \\ [N]=\eta}} \mathbf{P}^{\mathrm{f}}_{N,\omega}(S_{\lfloor tN \rfloor} > 0) = \left(1 - \frac{2 \arcsin \sqrt{t}}{\pi}\right) \mathrm{p}_\omega + \frac{2 \arcsin \sqrt{t}}{\pi} \mathrm{q}_{\omega,\eta}$$
(5.18)
$$=: \mathrm{p}^{\mathrm{f}}_{\omega,\eta}(t),$$

where $\mathsf{p}_\omega$ is the same as above [see (5.17)], while $\mathsf{q}_{\omega,\eta}$ is defined in (5.19) below. We again stress that we should actually fix the values of $G_{\lfloor tN \rfloor}$ and $D_{\lfloor tN \rfloor}$ as in (5.16), proving that the limiting probability is either $\mathsf{p}_\omega$ or $\mathsf{q}_{\omega,\eta}$ according to whether $D_{\lfloor tN \rfloor} \leq N$ or $D_{\lfloor tN \rfloor} > N$, but this will be clear from the steps below. Formula (5.18) follows from the fact that

$$\mathbf{P}^{\mathrm{f}}_{N,\omega}(S_{\lfloor tN \rfloor} > 0)$$
$$= \sum_{\alpha,\beta} \sum_{x<\lfloor tN \rfloor} \sum_{y>\lfloor tN \rfloor} \frac{\mathcal{Z}_{0,\alpha}(x)\rho^+_{\alpha,\beta}(y-x)M_{\alpha,\beta}(y-x)Z^{\mathrm{f}}_{N-y,\theta_{[y]}\omega}}{Z^{\mathrm{f}}_{N,\omega}}$$
$$+ \sum_\alpha \sum_{x<\lfloor tN \rfloor} \frac{\mathcal{Z}_{0,\alpha}(x)\rho^+_{\alpha,[N]}(N-x)P(N-x)\exp(\widetilde{\Phi}_{[x],[N]}(N-x))}{Z^{\mathrm{f}}_{N,\omega}}.$$

Letting $N \to \infty$ with $[N] = \eta$, by (3.27), the first term in the r.h.s. converges to

$$\int_0^t \frac{dx}{x^{1/2}} \int_t^1 \frac{dy}{(y-x)^{3/2}} \sum_{\alpha,\beta} \frac{1}{T^2} \frac{T^2 \xi_0 \zeta_\alpha}{2\pi \langle \zeta, L\xi \rangle}$$
$$\times c_K \frac{1}{2} \exp(\omega^{(0)}_\beta) \frac{\xi_\beta(T/2)\langle \zeta, \widetilde{L}_{\cdot,\eta} \rangle}{\langle \zeta, L\xi \rangle} \cdot \frac{\langle \zeta, L\xi \rangle}{\xi_0(T/2)\langle \zeta, \widetilde{L}_{\cdot,\eta} \rangle}$$
$$= \left(1 - \frac{2\arcsin\sqrt{t}}{\pi}\right) \cdot \mathsf{p}_\omega,$$

while the second term converges to

$$\int_0^t \frac{dx}{x^{1/2}(1-x)^{1/2}} \frac{1}{T} \sum_\alpha \frac{T^2 \xi_0 \zeta_\alpha}{2\pi \langle \zeta, L\xi \rangle} c_K \cdot \frac{\langle \zeta, L\xi \rangle}{\xi_0(T/2)\langle \zeta, \widetilde{L}_{\cdot,\eta} \rangle}$$
$$= \frac{2\arcsin\sqrt{t}}{\pi} \cdot \frac{c_K \sum_\gamma \zeta_\gamma}{\langle \zeta, \widetilde{L}_{\cdot,\eta} \rangle}.$$

Therefore, we obtain (5.18) with

(5.19) $$\mathsf{q}_{\omega,\eta} = \frac{c_K \sum_\gamma \zeta_\gamma}{\langle \zeta, \widetilde{L}_{\cdot,\eta} \rangle}.$$

REMARK 5.8. We observe that, by (3.24):

- if $h_\omega > 0$, or if $h_\omega = 0$ and $\Sigma \equiv 0$, then $\mathsf{p}^{\mathrm{f}}_{\omega,\eta}(t) = \mathsf{q}_{\omega,\eta} = \mathsf{p}_\omega$ for all $t$ and $\eta$;
- in the remaining case, that is, if $\omega \in \mathcal{P}$, in general, $\mathsf{p}^{\mathrm{f}}_{\omega,\eta}(t)$ depends on $t$ and $\eta$.

Now that we have proven the convergence of the probabilities of the signs of the excursion, in order to conclude the proof of point (2) of Theorem 1.3 and Proposition 2.1, it is enough to use the excursion theory of Brownian motion. For the details, we refer to the proof of Theorem 11 in [9]. □

## APPENDIX A: AN ASYMPTOTIC RESULT

In this appendix we will to prove that equation (3.15) holds true, but we first need some preliminary notation and results.

Given an *irreducible* $T \times T$ matrix $Q_{\alpha,\beta}$ with *nonnegative entries* [22], its Perron–Frobenius eigenvalue (= spectral radius) will be denoted by $\mathtt{Z} = \mathtt{Z}(Q)$ and the corresponding left and right eigenvectors (with any normalization) will be denoted by $\{\zeta_\alpha\}, \{\xi_\alpha\}$. We recall that $\zeta_\alpha, \xi_\alpha > 0$. Being a simple root of the characteristic polynomial, $\mathtt{Z}(Q)$ is an *analytic* function of the entries of $Q$, and

$$(\mathrm{A.1}) \qquad \frac{\partial \mathtt{Z}}{\partial Q_{\alpha,\beta}} = \frac{\zeta_\alpha \xi_\beta}{\langle \zeta, \xi \rangle}.$$

Hence, $\mathtt{Z}(Q)$ is a strictly increasing function of each of the entries of $Q$.

Now, let $Q$ denote the transition matrix of an irreducible, positive recurrent Markov chain and let us introduce the matrix $Q^{(\gamma)}$ and the vector $\delta^{(\gamma)}$, defined by

$$[Q^{(\gamma)}]_{\alpha,\beta} := Q_{\alpha,\beta} \mathbf{1}_{(\beta \neq \gamma)}, \qquad [\delta^{(\gamma)}]_\alpha := \mathbf{1}_{(\alpha = \gamma)}.$$

By monotonicity, $\mathtt{Z}(Q^{(\gamma)}) < \mathtt{Z}(Q) = 1$ for all $\gamma$. We can then define the geometric series

$$(1 - Q^{(\gamma)})^{-1} := \sum_{k=0}^{\infty} (Q^{(\gamma)})^k.$$

The interesting point is that, for every fixed $\gamma$, the vector $\alpha \mapsto [(1-Q^{(\gamma)})^{-1}]_{\gamma,\alpha}$ is (a multiple of) the left Perron–Frobenius eigenvector of the matrix $Q$. Similarly, the vector $\alpha \mapsto [(1 - Q^{(\gamma)})^{-1} \cdot Q]_{\alpha,\gamma}$ is (a multiple of) the right Perron–Frobenius eigenvector of $Q$. More precisely, we have

$$(\mathrm{A.2}) \qquad [(1 - Q^{(\gamma)})^{-1}]_{\gamma,\alpha} = \frac{\nu_\alpha}{\nu_\gamma}, \qquad [(1 - Q^{(\gamma)})^{-1} \cdot Q]_{\alpha,\gamma} = 1,$$

where $\{\nu_\alpha\}_\alpha$ is the unique invariant law of the chain, that is, $\sum_\alpha \nu_\alpha Q_{\alpha,\beta} = \nu_\beta$ and $\sum_\alpha \nu_\alpha = 1$. Equation (A.2) can be proven by exploiting its probabilistic interpretation in terms of the expected number of visits to state $\alpha$ before the first return to site $\gamma$; see [2], Section I.3.

Next we turn to our main problem. We recall, for convenience, the notation introduced in Sections 3.1 and 3.3. The process $\{\tau_k\}_{k \geq 0}$ where $\tau_0 = 0$ and $\tau_k = T_1 + \cdots + T_k$ is a Markov renewal process associated with the semi-Markov kernel $\Gamma_{\alpha,\gamma}(n)$ [defined in (3.13)] and $\{J_k\}_{k\geq 0}$ is its modulating chain. We denote by $\mathbb{P}_\beta$ the law of $\{(J_k, \tau_k)\}_{k\geq 0}$ with starting point $J_0 = \beta$ and we set $\ell := \inf\{k > 0 : J_k = \beta\}$. $q^\beta(x)$ then denotes the law of $\tau_\ell$ under $\mathbb{P}_\beta$ and we want to determine its asymptotic behavior.

We anticipate that the notation is necessarily quite involved, but the basic idea is simple. By the periodic structure of the kernel $\Gamma$, it follows that $q^\beta(x)$ is zero if $[x] \neq 0$. On the other hand, when $[x] = [0]$, by summing over the possible values of the index $\ell$ and using equation (3.6), we obtain

$$
\begin{aligned}
q^\beta(x) &= \mathbb{P}_\beta(\tau_1 = x, J_1 = \beta) \\
&\quad + \sum_{k=1}^\infty \mathbb{P}_\beta(J_i \neq \beta : 1 \leq i \leq k, J_{k+1} = \beta, \tau_{k+1} = x) \\
&= \sum_{k=0}^\infty ((\Gamma^{(\beta)})^{*k} * \Gamma)_{\beta,\beta}(x),
\end{aligned}
\tag{A.3}
$$

where we have introduced the kernel $\Gamma_{\alpha,\gamma}^{(\beta)}(x) := \Gamma_{\alpha,\gamma}(x)\mathbf{1}_{(\gamma \neq \beta)}$ that gives the law of the steps with index $k < \ell$. Looking at (A.3), we set $V_{\alpha,\gamma}^{(\beta)}(x) := \sum_{k=0}^\infty ((\Gamma^{(\beta)})^{*k})_{\alpha,\gamma}(x)$ and we can write

$$
q^\beta(x) = (V^{(\beta)} * \Gamma)_{\beta,\beta}(x) = \sum_{\gamma \in \mathbb{S}} \sum_{y=0}^{x-1} V_{\beta,\gamma}^{(\beta)}(y) \Gamma_{\gamma,\beta}(x-y).
\tag{A.4}
$$

The asymptotic behavior of $q^\beta(x)$ can be extracted from the above expression. To this end, we need to know both the asymptotic behavior as $n \to \infty$ and the sum over $n \in \mathbb{N}$ of the two kernels $\Gamma_{\gamma,\beta}(n)$ and $V_{\beta,\gamma}^{(\beta)}(n)$ appearing in the r.h.s.

- By (3.13) and (3.10), as $n \to \infty$ along $[n] = \beta - \gamma$ we have

$$
\Gamma_{\gamma,\beta}(n) \sim \frac{\widehat{L}_{\gamma,\beta}}{n^{3/2}}, \qquad \text{where } \widehat{L}_{\gamma,\beta} := L_{\gamma,\beta} \frac{\xi_\beta}{\xi_\gamma}.
\tag{A.5}
$$

Moreover, the sum over $n \in \mathbb{N}$ gives

$$
\sum_{n \in \mathbb{N}} \Gamma_{\gamma,\beta}(n) = B_{\gamma,\beta} \frac{\xi_\beta}{\xi_\gamma} =: \widehat{B}_{\gamma,\beta}.
\tag{A.6}
$$

- For the asymptotic behavior of the kernel $V^{(\beta)} := \sum_{k=0}^\infty (\Gamma^{(\beta)})^{*k}$, we can apply the theory developed in Section 3.4 for the case $\delta^\omega < 1$ because the matrix $\sum_{x \in \mathbb{N}} \Gamma_{\alpha,\gamma}^{(\beta)}(x)$ is just $[\widehat{B}^{(\beta)}]_{\alpha,\gamma}$ by (A.6) (we recall the convention $[Q^{(\beta)}]_{\alpha,\gamma} := Q_{\alpha,\gamma} \mathbf{1}_{(\gamma \neq \beta)}$ for any matrix $Q$) which has Perron–Frobenius eigenvalue strictly smaller than 1. Since

$$\Gamma^{(\beta)}_{\alpha,\gamma}(n) \sim \frac{[\widehat{L}^{(\beta)}]_{\alpha,\gamma}}{n^{3/2}}, \qquad n \to \infty, [n] = \gamma - \alpha,$$

we can apply (3.18) to obtain the asymptotic behavior as $n \to \infty$, $[n] = \alpha - \gamma$:

(A.7) $$V^{(\beta)}_{\alpha,\gamma}(n) \sim ([(1-\widehat{B}^{(\beta)})^{-1}\widehat{L}^{(\beta)}(1-\widehat{B}^{(\beta)})^{-1}]_{\alpha,\gamma})\frac{1}{n^{3/2}}.$$

On the other hand, for the sum over $n \in \mathbb{N}$, an analog of (3.19) yields

(A.8) $$\sum_{n \in \mathbb{N}} V^{(\beta)}_{\alpha,\gamma}(n) = \sum_{k=0}^{\infty} [(\widehat{B}^{(\beta)})^k]_{\alpha,\gamma} = [(1-\widehat{B}^{(\beta)})^{-1}]_{\alpha,\gamma}.$$

As equations (A.5) and (A.7) show, both kernels $V^{(\beta)}$ and $\Gamma$ have an $n^{-3/2}$ tail. From (A.4), it then follows that as $x \to \infty$ along $[x] = 0$,

$$q^{\beta}(x) \sim \sum_{\gamma \in \mathbb{S}} \left\{ \left(\sum_{n \in \mathbb{N}} V^{(\beta)}_{\beta,\gamma}(n)\right) \Gamma_{\gamma,\beta}(x) + V^{(\beta)}_{\beta,\gamma}(x) \left(\sum_{n \in \mathbb{N}} \Gamma_{\gamma,\beta}(n)\right) \right\}.$$

It now suffices to apply (A.8), (A.5), (A.7) and (A.6) to see that, indeed, $q^{\beta}(x) \sim c_{\beta}/x^{3/2}$ as $x \to \infty$ along $[x] = 0$, where the positive constant $c_{\beta}$ is given by

$$c_{\beta} = [(1-\widehat{B}^{(\beta)})^{-1} \cdot \widehat{L}]_{\beta,\beta} + [(1-\widehat{B}^{(\beta)})^{-1} \cdot \widehat{L}^{(\beta)} \cdot (1-\widehat{B}^{(\beta)})^{-1} \cdot \widehat{B}]_{\beta,\beta}.$$

Using the fact that $[(1-\widehat{B}^{(\beta)})^{-1} \cdot \widehat{B}]_{\beta,\beta} = 1$, which follows from (A.2) applied to the matrix $Q = \widehat{B}$, we can rewrite the above expression as

$$c_{\beta} = [(1-\widehat{B}^{(\beta)})^{-1} \cdot \widehat{L} \cdot (1-\widehat{B}^{(\beta)})^{-1} \cdot \widehat{B}]_{\beta,\beta} = \frac{1}{\nu_{\beta}} \sum_{\alpha,\gamma \in \mathbb{S}} \nu_{\alpha} \widehat{L}_{\alpha,\gamma},$$

where $\{\nu_{\alpha}\}_{\alpha}$ is the invariant measure of the matrix $\widehat{B}$ and the second equality again follows from (A.2). However, from (A.6) it is easily seen that $\{\nu_{\alpha}\} = \{\zeta_{\alpha}\xi_{\alpha}\}$ and recalling the definition (A.5) of $\widehat{L}$, we finally obtain the expression for $c_{\beta}$ given in equation (3.15):

(A.9) $$c_{\beta} = \frac{1}{\zeta_{\beta}\xi_{\beta}} \sum_{\alpha,\gamma} \zeta_{\alpha} L_{\alpha,\gamma} \xi_{\gamma}.$$

## APPENDIX B: A LOCALIZATION ARGUMENT

Let us give a proof that for the *copolymer near a selective interface* model, described in Section 1.1, the charge $\omega$ never belongs to $\mathcal{P}$ [see (2.9) for the definition of $\mathcal{P}$]. More precisely, we will show that if $h_\omega = 0$ and $\Sigma \not\equiv 0$, then $\delta^\omega > 1$. That is, the periodic copolymer with zero-mean, nontrivial charges is always localized. As a matter of fact, this is an immediate consequence of the estimates on the critical line obtained in [5]. However, we want to give here an explicit proof, both because it is more direct and because the model studied in [5] is built over the simple random walk measure, corresponding to $p = 1/2$ with the language of Section 1, while we consider the case $p < 1/2$.

We recall that, by (A.1), the Perron–Frobenius eigenvalue $\mathtt{Z}(Q)$ of an irreducible matrix $Q$ is increasing in the entries of $Q$. We also point out a result proved by Kingman [19]: if the matrix $Q = Q(t)$ is a function of a real parameter $t$ such that all the entries $Q_{\alpha,\beta}(t)$ are *log-convex* functions of $t$ [i.e., $t \mapsto \log Q_{\alpha,\beta}(t)$ is convex for all $\alpha, \beta$], then $t \mapsto \mathtt{Z}(Q(t))$ is also a log-convex function of $t$.

Next, we come to the *copolymer near a selective interface* model: with reference to the general Hamiltonian (1.3), we are assuming that $\omega_n^{(0)} = \widetilde{\omega}_n^{(0)} = 0$ and $h_\omega = 0$ [where $h_\omega$ was defined in (1.5)]. In this case, the integrated Hamiltonian $\Phi_{\alpha,\beta}(\ell)$ [see (2.4)] is given by

$$\Phi_{\alpha,\beta}(\ell) = \begin{cases} 0, & \text{if } \ell = 1 \text{ or } \ell \notin \beta - \alpha, \\ \log[\frac{1}{2}(1 + \exp(\Sigma_{\alpha,\beta}))], & \text{if } \ell > 1 \text{ and } \ell \in \beta - \alpha. \end{cases}$$

We recall that the law of the first return to zero of the original walk is denoted by $K(\cdot)$ [see (2.1)] and we introduce the function $q : \mathbb{S} \to \mathbb{R}^+$ defined by

$$q(\gamma) := \sum_{x \in \mathbb{N}, [x] = \gamma} K(x)$$

[notice that $\sum_\gamma q(\gamma) = 1$]. The matrix $B_{\alpha,\beta}$ defined by (2.6) then becomes

$$(\text{B.1}) \qquad B_{\alpha,\beta} = \begin{cases} \frac{1}{2}(1 + \exp(\Sigma_{\alpha,\beta}))q(\beta - \alpha), \\ \quad \text{if } \beta - \alpha \neq [1], \\ K(1) + \frac{1}{2}(1 + \exp(\Sigma_{\alpha,\alpha+[1]})) \cdot (q([1]) - K(1)), \\ \quad \text{if } \beta - \alpha = [1]. \end{cases}$$

By (2.7), to prove localization, we must show that the Perron–Frobenius eigenvalue of the matrix $(B_{\alpha,\beta})$ is strictly greater than 1, that is, $\mathtt{Z}(B) > 1$. Applying the elementary convexity inequality $(1 + \exp(x))/2 \geq \exp(x/2)$ to (B.1), we get

$$(\text{B.2}) \qquad B_{\alpha,\beta} \geq \widetilde{B}_{\alpha,\beta} := \begin{cases} \exp(\Sigma_{\alpha,\beta}/2)q(\beta - \alpha) \\ \quad \text{if } \beta - \alpha \neq [1], \\ K(1) + \exp(\Sigma_{\alpha,\alpha+[1]}/2) \cdot (q([1]) - K(1)) \\ \quad \text{if } \beta - \alpha = [1]. \end{cases}$$

By hypothesis, $\Sigma_{\alpha_0,\beta_0} \neq 0$ for some $\alpha_0, \beta_0$, therefore the inequality above is strict for $\alpha = \alpha_0$, $\beta = \beta_0$. We have already observed that the Perron–Frobenius eigenvalue is a strictly increasing function of the entries of the matrix, hence $\mathsf{Z}(B) > \mathsf{Z}(\widetilde{B})$. Therefore, it only remains to show that $\mathsf{Z}(\widetilde{B}) \geq 1$ and the proof will be complete.

Again, an elementary convexity inequality applied to the second line of (B.2) yields

(B.3) $$\widetilde{B}_{\alpha,\beta} \geq \widehat{B}_{\alpha,\beta} := \exp(c(\beta - \alpha)\Sigma_{\alpha,\beta}/2) \cdot q(\beta - \alpha),$$

where

$$c(\gamma) := \begin{cases} 1, & \text{if } \gamma \neq [1], \\ \dfrac{q([1]) - K(1)}{q([1])}, & \text{if } \gamma = [1]. \end{cases}$$

We will prove that $\mathsf{Z}(\widehat{B}) \geq 1$. Observe that by setting $v_\alpha := \Sigma_{[0],\alpha}$, we can write

$$\Sigma_{\alpha,\beta} = \Sigma_{[0],\beta} - \Sigma_{[0],\alpha} = v_\beta - v_\alpha.$$

We then make a similarity transformation via the matrix $L_{\alpha,\beta} := \exp(v_\beta/2)\mathbf{1}_{(\beta=\alpha)}$, getting

$$C_{\alpha,\beta} := [L \cdot \widehat{B} \cdot L^{-1}]_{\alpha,\beta} = \exp((c(\beta - \alpha) - 1)\Sigma_{\alpha,\beta}/2) \cdot q(\beta - \alpha)$$
$$= \exp(d\Sigma_{\alpha,\alpha+[1]}\mathbf{1}_{(\beta-\alpha=1)}) \cdot q(\beta - \alpha),$$

where we have introduced the constant $d := -K(1)/(2q([1]))$. Of course, $\mathsf{Z}(\widehat{B}) = \mathsf{Z}(C)$. Also, notice that by the very definition of $\Sigma_{\alpha,\beta}$, we have $\Sigma_{\alpha,\alpha+[1]} = \omega^{(-1)}_{\alpha+[1]} - \omega^{(+1)}_{\alpha+[1]}$, hence the hypothesis $h_\omega = 0$ yields $\sum_{\alpha \in \mathbb{S}}(\Sigma_{\alpha,\alpha+[1]}) = 0$.

Thus we are finally left with the task of showing that $\mathsf{Z}(C) \geq 1$, where $C_{\alpha,\beta}$ is an $\mathbb{S} \times \mathbb{S}$ matrix of the form

$$C_{\alpha,\beta} = \exp(w_\alpha \mathbf{1}_{(\beta-\alpha=1)}) \cdot q(\beta - \alpha), \qquad \text{where } \sum_\alpha w_\alpha = 0, \sum_\gamma q(\gamma) = 1.$$

To this end, we introduce an interpolation matrix

$$C_{\alpha,\beta}(t) := \exp(t \cdot w_\alpha \mathbf{1}_{(\beta-\alpha=1)}) \cdot q(\beta - \alpha),$$

defined for $t \in \mathbb{R}$, and notice that $C(1) = C$. Let us denote by $\eta(t) := \mathsf{Z}(C(t))$ the Perron–Frobenius eigenvalue of $C(t)$. As the entries of $C(t)$ are log-convex functions of $t$, it follows that $\eta(t)$ is also log-convex, therefore in particular, convex. Moreover, $\eta(0) = 1$ [the matrix $C(0)$ is bistochastic] and using (A.1), one easily checks that $\frac{\mathrm{d}}{\mathrm{d}t}\eta(t)|_{t=0} = 0$. Since, clearly, $\eta(t) \geq 0$ for all $t \in \mathbb{R}$, by convexity, it follows that, indeed, $\eta(t) \geq 1$ for all $t \in \mathbb{R}$ and the proof is complete.

# REFERENCES

[1] ALEXANDER, K. S. and SIDORAVICIUS, V. (2006). Pinning of polymers and interfaces by random potentials. *Ann. Appl. Probab.* **16** 636–669. MR2244428
[2] ASMUSSEN, S. (2003). *Applied Probability and Queues*, 2nd ed. Springer, New York. MR1978607
[3] BOLTHAUSEN, E. (1976). On a functional central limit theorem for random walks conditioned to stay positive. *Ann. Probab.* **4** 480–485. MR0415702
[4] BOLTHAUSEN, E. and DEN HOLLANDER, F. (1997). Localization transition for a polymer near an interface. *Ann. Probab.* **25** 1334–1366. MR1457622
[5] BOLTHAUSEN, E. and GIACOMIN, G. (2005). Periodic copolymers at selective interfaces: A large deviations approach. *Ann. Appl. Probab.* **15** 963–983. MR2114996
[6] CARAVENNA, F., GIACOMIN, G. and ZAMBOTTI, L. (2006). Sharp asymptotic behavior for wetting models in $(1+1)$-dimension. *Electron. J. Probab.* **11** 345–362. MR2217821
[7] CARAVENNA, F., GIACOMIN, G. and ZAMBOTTI, L. (2007). Infinite volume limits of polymer chains with periodic charges. *Markov Process. Related Fields.* To appear. Available at www.arxiv.org/abs/math.PR/0604426.
[8] CARAVENNA, F., GIACOMIN, G. and ZAMBOTTI, L. (2007). Tightness conditions for polymer measures. Available at www.arxiv.org/abs/math.PR/0702331.
[9] DEUSCHEL, J.-D. GIACOMIN, G. and ZAMBOTTI, L. (2005). Scaling limits of equilibrium wetting models in $(1+1)$-dimension. *Probab. Theory Related Fields* **132** 471–500. MR2198199
[10] DONEY, R. A. (1997). One-sided local large deviation and renewal theorems in the case of infinite mean. *Probab. Theory Related Fields* **107** 451–465. MR1440141
[11] FELLER, W. (1971). *An Introduction to Probability Theory and Its Applications* **II**, 2nd ed. Wiley, New York. MR0270403
[12] FITZSIMMONS, P. J., FRISTEDT, B. and MAISONNEUVE, B. (1985). Intersections and limits of regenerative sets. *Z. Wahrsch. Verw. Gebiete* **70** 157–173. MR0799144
[13] GALLUCCIO, S. and GRABER, R. (1996). Depinning transition of a directed polymer by a periodic potential: A $d$-dimensional solution. *Phys. Rev. E* **53** R5584–R5587.
[14] GARSIA, A. and LAMPERTI, J. (1963). A discrete renewal theorem with infinite mean. *Comment. Math. Helv.* **37** 221–234. MR0148121
[15] GIACOMIN, G. (2007). *Random Polymer Models*. Imperial College Press, London.
[16] GIACOMIN, G. and TONINELLI, F. L. (2005). Estimates on path delocalization for copolymers at selective interfaces. *Probab. Theory Related Fields* **133** 464–482. MR2197110
[17] ISOZAKI, Y. and YOSHIDA, N. (2001). Weakly pinned random walk on the wall: Pathwise descriptions of the phase transition. *Stochastic Process. Appl.* **96** 261–284. MR1865758
[18] KAIGH, W. D. (1976). An invariance principle for random walk conditioned by a late return to zero. *Ann. Probab.* **4** 115–121. MR0415706
[19] KINGMAN, J. F. C. (1961). A convexity property of positive matrices. *Quart. J. Math. Oxford Ser.* (*2*) **12** 283–284. MR0138632
[20] LE GALL, J.-F. and ROSEN, J. (1991). The range of stable random walks. *Ann. Probab.* **19** 650–705. MR1106281
[21] MATHERON, G. (1975). *Random Sets and Integral Geometry*. Wiley, New York. MR0385969
[22] MINC, H. (1988). *Nonnegative Matrices*. Wiley, New York. MR0932967

[23] MONTHUS, C., GAREL, T. and ORLAND, H. (2000). Copolymer at a selective interface and two dimensional wetting: A grand canonical approach. *Eur. Phys. J. B* **17** 121–130.
[24] NAIDEDOV, A. and NECHAEV, S. (2001). Adsorption of a random heteropolymer at a potential well revisited: Location of transition point and design of sequences. *J. Phys. A*: *Math. Gen.* **34** 5625–5634. MR1857165
[25] NECHAEV, S. and ZHANG, Y.-C. (1995). Exact solution of the 2D wetting problem in a periodic potential. *Phys. Rev. Lett.* **74** 1815–1818.
[26] REVUZ, D. and YOR, M. (1991). *Continuous Martingales and Brownian Motion*. Springer, Berlin. MR1083357
[27] SOMMER, J.-U. and DAOUD, M. (1995). Copolymers at selective interfaces. *Europhys. Lett.* **32** 407–412.
[28] SOTEROS, C. E. and WHITTINGTON, S. G. (2004). The statistical mechanics of random copolymers. *J. Phys. A*: *Math. Gen.* **37** R279–R325. MR2097625

F. CARAVENNA
DIPARTIMENTO DI MATEMATICA PURA E APPLICATA
UNIVERSITÀ DEGLI STUDI DI PADOVA
VIA TRIESTE 63
35121 PADOVA
ITALY
E-MAIL: francesco.caravenna@math.unipd.it

G. GIACOMIN
LABORATOIRE DE PROBABILITÉS
ET MODÈLES ALÉATOIRES (CNRS U.M.R. 7599)
AND
UNIVERSITÉ PARIS 7—DENIS DIDEROT
U.F.R. MATHEMATIQUES
CASE 7012
2 PLACE JUSSIEU
75251 PARIS CEDEX 05
FRANCE
E-MAIL: giacomin@math.jussieu.fr

L. ZAMBOTTI
LABORATOIRE DE PROBABILITÉS
ET MODÈLES ALÉATOIRES (CNRS U.M.R. 7599)
AND
UNIVERSITÉ PARIS 6—PIERRE ET MARIE CURIE
U.F.R. MATHEMATIQUES
CASE 188
4 PLACE JUSSIEU
75252 PARIS CEDEX 05
FRANCE
E-MAIL: zambotti@ccr.jussieu.fr